\documentclass[12p]{article}
\usepackage{amsmath, amsthm}
\usepackage{amsfonts} %,showkeys}
\usepackage{amssymb}
\usepackage{graphicx}
\usepackage[monochrome]{color}
\usepackage{hyperref}

\newtheorem{thm}{Theorem}[section]
\newtheorem{cor}[thm]{Corollary}
\newtheorem{lem}[thm]{Lemma}
\newtheorem{prop}[thm]{Proposition}
\newtheorem{defn}[thm]{Definition}

\numberwithin{equation}{section}

\newcommand{\dx}{\,{\rm d}x}
\newcommand{\dy}{\,{\rm d}y}

\newcommand{\dt}{\,{\rm d}t}
\newcommand{\rd}{{\rm d}}

\newcommand{\dom}{\mathrm{dom}}
\newcommand{\A}{\mathcal L}

\newcommand{\AM}{\mathcal{L}^{\frac{1}{2}}}
\newcommand{\AI}{\mathcal{L}^{-1}}
\newcommand{\AIM}{\mathcal{L}^{-\frac{1}{2}}}

\renewcommand{\k}{\kappa}
\newcommand{\K}{{\mathbb K}}

\def\quotient#1#2{%
    \raise1ex\hbox{$#1$}\Big/\lower1ex\hbox{$#2$}%
}

\def\LL{\mathrm{L}} %per gli spazi L^p
\def\supp{\mathrm{supp}} %per il supporto
\newcommand{\RR}{\mathbb{R}}
\newcommand{\NN}{\mathbb{N}}

\newcommand{\w}{\overline{w}}
\newcommand{\tr}{{\rm Tr}}

\def\ee{\mathrm{e}} %per l'esponenziale
\def\dist{\mathrm{dist}} %per la distanza
\def\qed{\,\unskip\kern 6pt \penalty 500
\raise -2pt\hbox{\vrule \vbox to8pt{\hrule width 6pt
\vfill\hrule}\vrule}\par}
\definecolor{darkblue}{rgb}{0.05, .05, .65}
\definecolor{darkgreen}{rgb}{0.05, .70, .05}
\definecolor{darkred}{rgb}{0.8,0,0}
\begin{document}
\title{ \bf Existence, Uniqueness and Asymptotic behaviour\\ for fractional porous medium equations\\ on bounded domains }
\author{Matteo Bonforte$^{\,a,\,b}$, %\footnote{e-mail address:~bonforte@calvino.polito.it},
Yannick Sire$^{\,c,\,d}$,
~and~ Juan Luis V\'azquez$^{\,a,\,e}$} %\footnote{e-mail address:~juanluis.vazquez@uam.es}}
\date{} %%  this cancels date in article format

\maketitle

\begin{abstract}
 We consider nonlinear diffusive evolution equations posed on bounded space domains, governed by  fractional Laplace-type operators, and involving porous medium type nonlinearities. We establish  existence and uniqueness results in a suitable class of solutions using the theory of maximal monotone operators on dual spaces. Then we describe the long-time asymptotics in terms of separate-variables solutions of the friendly giant type. As a by-product, we obtain an existence and uniqueness result for semilinear elliptic non local equations with sub-linear nonlinearities.  The Appendix contains a review of the theory of fractional Sobolev spaces and of the interpolation theory that are used in the rest of the paper.

\end{abstract}

\vfill



\noindent (a) Departamento de Matem\'{a}ticas, Universidad
Aut\'{o}noma de Madrid, Campus de Cantoblanco, 28049 Madrid, Spain\\
\noindent (b) e-mail address:~matteo.bonforte@uam.es \\
\noindent (c) Universit\'e Aix-Marseille,  I2M, Centre de Math\'ematique et Informatique, Technop\^ole de Chateau-Giombert, Marseille, France \\
\noindent (d) e-mail address:~sire@cmi.univ-mrs.fr\\
\noindent (e) e-mail address:~juanluis.vazquez@uam.es\\
% 35B (1973-now) Qualitative properties of solutions
% 35B45 (1973-now) A priori estimates
% 35B65 (1980-now) Smoothness and regularity of solutions of PDE
% 35K (1973-now) Parabolic equations and systems
%35K55 (1973-now) Nonlinear PDE of parabolic type
%35K65 (1980-now) Parabolic partial differential equations of degenerate type
%   35K61   	Nonlinear initial-boundary value problems for nonlinear parabolic equations
%   35B40   	Asymptotic behavior of solutions

% 	35A01   	Existence problems: global existence, local existence, non-existence
%	35A02   	Uniqueness problems: global uniqueness, local uniqueness, non-uniqueness

\newpage
\small
\tableofcontents

\newpage

\normalsize
\section{Introduction}

Nonlocal equations have attracted much attention in the last decade. The basic operator involved is the so-called fractional Laplacian $(-\Delta)^s$, $s \in (0,1)$ in $\mathbb R^d$, a Fourier multiplier with symbol $|\xi|^{2s}$. This operator and a number of variants appear in numerous areas of mathematics and mathematical physics: harmonic analysis, probability theory, potential theory, quantum mechanics, statistical physics, cosmology,...
From the point of view of mathematical analysis, a wide range of elliptic and parabolic problems involving such operators has been studied.

We are interested in the study of nonlinear diffusion processes  involving fractional Laplacian operators for which there is a growing literature, cf. \cite{Vaz2012abel}, \cite{Vaz2014}. In this paper we focus our attention on problems posed on bounded space domains. More precisely, we consider the Fractional Diffusion Equation (FDE)
\begin{equation}\label{FPME.equation}
\partial_t u +\A(u^m)=0  \qquad\mbox{in }(0,\infty)\times \Omega\,,
\end{equation}
where $\Omega\subset \RR^d$ is a bounded domain with smooth boundary, $m>1$, $0<s<1$ and $d\ge 1$. The linear operator $\A$ is a fractional power of the Laplacian subject to suitable Dirichlet boundary conditions. We point out that the proper definition of the fractional operator $\A$  on a bounded space domain is not immediate and offers some choices. Actually, we consider  two essentially different operators corresponding to zero boundary conditions, and present a unified theory  which works for both operators modulo appropriately choosing  the functional setting, which turns out to be a delicate part of the work. Therefore, we devote two long final appendices to clarify the role of the fractional Sobolev spaces and the interpolation theory that is required.

Using these tools, we  establish existence, uniqueness and different estimates for a suitable class of solutions of \eqref{FPME.equation} with appropriate initial conditions (and zero boundary conditions as mentioned). These solutions exist globally in time. Let us mention that our method allows to prove existence and uniqueness for a wider class of linear operators and nonlinearities, and give examples of such extensions.

In a second stage, we describe in detail the large-time behaviour of such solutions, with convergence to equilibrium after a natural time-dependent rescaling.

As a by-product of our method, we will obtain an existence and uniqueness result for nonlocal elliptic equations of the type
\begin{equation}
\A v =v^p\qquad\mbox{in } \Omega
\end{equation}
where $p=1/m<1$, i.e. a sublinear growth. This topic has an independent interest; for instance we obtain in this way a unique stationary state for the evolution equation
$v_t+ \A v =v^p$ with $0<p<1.$

%%%%%%%%%%%%%%%%%%%%%%%%%%%%%%%%%%%%%%%%%%%%%%%%%%%%%%%%%%%%%%%%%%%%%%%%%%%%%%%%%%

\section{Operators, evolution problem, and main results}\label{sec.main}

\subsection{The operators}

\noindent $\bullet$ {\sl The spectral Laplacian: } If one considers the classical Dirichlet Laplacian $\Delta_{\Omega}$ on the domain $\Omega$\,, then the {spectral definition} of the fractional power of $\Delta_{\Omega}$ relies on the following formulas:
\begin{equation}\label{sLapl.Omega.Spectral}
\displaystyle(-\Delta_{\Omega})^{s}
g(x)=\sum_{j=1}^{\infty}\lambda_j^s\, \hat{g}_j\, \phi_j(x)
=\frac1{\Gamma(-s)}\int_0^\infty
\left(e^{t\Delta_{\Omega}}g(x)-g(x)\right)\frac{dt}{t^{1+s}}.
\end{equation}
Here $\lambda_j>0$, $j=1,2,\ldots$ are the eigenvalues of the Dirichlet Laplacian on $\Omega$ with zero boundary conditions\,, written in increasing order and repeated according to their multiplicity and $\phi_j$ are the corresponding normalized eigenfunctions, namely
\[
\hat{g}_j=\int_\Omega g(x)\phi_j(x)\dx\,,\qquad\mbox{with}\qquad \|\phi_j\|_{\LL^2(\Omega)}=1\,.
\]
The first part of the formula is therefore an interpolation definition. The second part gives an equivalent definition in terms of the semigroup associated to the Laplacian.
We will denote the operator defined in such a way as $\A_{1,s}=(-\Delta_{\Omega})^s$\,, and call it the \textit{spectral fractional Laplacian}. The reader should remember that the zero  boundary conditions are built into the definition of the operator. There is another way of defining the SFL using the Caffarelli-Silvestre extension, which turns out to be equivalent, see e.g. \cite{CabTan, CaSi, c-d-d-s, Tan}\,.

We refer to \cite{DPQRV1,DPQRV2} for the basic theory of weak solutions to the initial and boundary value problem for equation \eqref{FPME.equation} with $u_0\in\LL^1(\Omega)$, as well as for information on the physical motivation and relevance of this nonlocal model.

\medskip

\noindent $\bullet$~{\sl The restricted fractional laplacian:} On the other hand, one can define a fractional Laplacian operator by using the integral representation in terms of hypersingular kernels
\begin{equation}\label{sLapl.Rd.Kernel}
(-\Delta_{\RR^d})^{s}  g(x)= c_{d,s}\mbox{
P.V.}\int_{\mathbb{R}^d} \frac{g(x)-g(z)}{|x-z|^{d+2s}}\,dz,
\end{equation}
where $c_{d,s}>0$ is a normalization constant, see \cite{DPQRV2}. In this case we materialize the zero Dirichlet condition  by restricting the operator to act only on functions that are zero outside $\Omega$. We will call  the operator defined in such a way the \textit{restricted fractional Laplacian} and use the specific notation $\A_{2,s}=(-\Delta_{|\Omega})^s$ when needed. In this case, the initial and boundary conditions associated to the fractional diffusion equation \eqref{FPME.equation} read
\begin{equation}\label{FPME.Dirichlet.conditions.Restricted}
\left\{
\begin{array}{lll}
u(t,x)=0\,,\; &\mbox{in }(0,\infty)\times\RR^d\setminus \Omega\,,\\
u(0,\cdot)=u_0\,,\; &\mbox{in }\Omega\,.
\end{array}
\right.
\end{equation}
So defined,  $\A_{2,s}$ is a self-adjoint operator on $\LL^2(\Omega)$\,, with a discrete spectrum: we will denote by $\lambda_{s, j}>0$, $j=1,2,\ldots$ its eigenvalues written in increasing order and repeated according to their multiplicity and we will denote by $\{\phi_{s, j}\}_j$ the corresponding set of eigenfunctions, normalized in $L^2(\Omega)$. In this case, the eigenvalues are  smaller  than the ones of the spectral Laplacian cf. \cite{ChSo}, and the corresponding eigenfunctions are known to be H\"older  only    continuous up to the boundary cf. \cite{Grub1, Grub2, RosSer, RosSer1, SV1}.

\noindent $\bullet$~{\sl Common notation.}
In the sequel we use $\A$ to refer to any of the two types of operators $\A_{1,s}$ or $\A_{2,s}$, $0<s<1$. Each one is defined on a Hilbert space
\begin{equation} \label{defH}
H(\Omega)=\{u=\sum_{k=1}^\infty u_k \phi_{s,k} \in L^2(\Omega)\; : \; \| u \|^2_{H} = \sum_{k=1}^\infty \lambda_{s,k}\vert u_k\vert^2 <+\infty\}\subset L^2(\Omega)
\end{equation}
with values in its dual $H^*$. The notation in the formula copies the one just used for the second operator. When applied to the first one we put here  $\phi_{s,k}=\phi_k$, and $\lambda_{s,k}=\lambda_{k}^s$.
In  Section \ref{sec.exist} we will review the functional theory and identify in a clear way what is exactly the space $H$, since we will need to prove some compactness results to have existence of solutions.
 Note that $H(\Omega)$ depends in principle on the type of operator and on the exponent $s$. It turns out that $H_{1,s}(\Omega)=H_{2,s}(\Omega)$ for each $s$, see Section \ref{sec.exist} below. We also remark that $H^*$ can be described as
the completion of the finite sums of the form
\begin{equation}
\label{defH*}
f=\sum_{k=1}^N c_k \phi_{s,k}
\end{equation}
with respect to the dual norm
\begin{equation}
 \| f \|^2_{H^*} = \sum_{k=1}^\infty \lambda_{s,k}^{-1}\vert c_k\vert^2\,,
\end{equation}
and it is a space of distributions. Moreover, the operator $\A$ is an isomorphism between $H$ and $H^*$, given by its action on the eigen-functions. If $u,v\in H$ and $f={\cal L}u$ we have, after this isomorphism,
$$
\langle f,v\rangle_{H^*\times H}=\langle u,v\rangle_{H\times H}=\sum_1^{\infty}\lambda_{s,k} u_kv_k.
$$
If it also happens that $f\in L^2(\Omega)$, then clearly we get $\langle f,v\rangle_{H^*\times H}=\int_{\Omega}fv\,\dx.$

\subsection{Evolution Problem. Existence and uniqueness of solutions}

In this paper we will consider the homogeneous Dirichlet problem for equation \eqref{FPME.equation}
where $\Omega\subset \RR^d$ is a bounded domain with smooth boundary, $m>1$, $0<s<1$ and $d\ge 1$\,.
For the general theory we will use signed solutions $u(t,x)$. In that case we will still use the notation $u^m=|u|^{m-1}u$ for the sake of brevity.

The definition of solution to be introduced below includes zero Dirichlet boundary conditions which are implicitly  taken as a consequence of the definition of both operators $\A$. We also need initial conditions
at $t=0$.

Our first goal is to prove an existence and uniqueness result. Actually, the result holds for more general nonlinearities than the power one at almost no extra cost so that we will adopt the more general context and consider the problem
\begin{equation}\label{FPME.Problem.PHI}
\left\{\begin{array}{lll}
u_t+\A(\varphi(u))=0 &  ~ {\rm in}~ (0,+\infty)\times \Omega\\
u(0,x)=u_0(x) & ~{\rm in}~ \Omega \\
% u(t,x)=0 & ~\mbox{ on $(0,+\infty)\times \Gamma$}
\end{array}
\right.
\end{equation}
where $\varphi:\RR\to\RR$ is a continuous, smooth and increasing function. We assume moreover that $\varphi'>0$, $\varphi(\pm\infty)=\pm\infty$ and $\varphi(0)=0$. The leading example will be $\varphi(u)=|u|^{m-1}u$ with $m>0$\,. We recall that zero Dirichlet conditions are built into the definition of the operator, and they have different meanings in the two cases we consider.

\medskip

Choosing the correct type of generalized solutions will be critical in the proofs. Here we use the concept of $H^*$ solution.

\begin{defn}\label{def.H*sols}
$u \in C([0,T], H^*(\Omega))$ is an $H^*$-solution if $\varphi(u) \in L^1([0,T],  H(\Omega) )$ such that
\begin{equation}\label{testH*}
\int_0^T \int_\Omega u\psi_t= \int_0^T \int_\Omega \varphi(u) \A \psi\,\qquad\forall\,\psi \in C^1_c([0,T], H(\Omega))
\end{equation}
\end{defn}

Notice that since $\A$ is an isomorphism from $H$ into $H^*$, then equation \eqref{testH*} is equivalent to $\psi =(\A)^{-1} \tilde \psi$

\begin{equation}
\int_0^T \int_\Omega u \partial_t ((\A)^{-1} \tilde \psi)= \int_0^T \int_\Omega \varphi(u) \tilde \psi\,\qquad\forall\,\tilde \psi \in C^1_c([0,T], H^*(\Omega))
\end{equation}
 This is a formulation of weak solutions for the potential equation $\partial_t \AI u+\varphi(u)=0$\,.

The main result on existence and uniqueness is as follows:

\begin{thm}\label{cor3.2}
For every $u_0\in H^*(\Omega)$ there exists a unique solution $u\in C([0,T]:H^*(\Omega))$ of  Problem {\rm \ref{FPME.Problem.PHI}} for every $T>0$, i.e. the solution is global in time. We also have
\begin{equation}\label{strong}
t\,\varphi(u)\in\LL^{\infty}(0,T\,:\,H^*(\Omega)),\quad
t\,\partial_tu\in\LL^{\infty}(0,T\,:\,H^*(\Omega)).
\end{equation}
We also  have $u\varphi(u)\in L^{1}((0,T)\times\Omega)$. The solution map $S_t: u_0\mapsto u(t)$ defines a semigroup of (non-strict) contractions in $H^*(\Omega)$, i.\,e.,
\begin{equation}\label{contractivity.H*}
\|u(t)-v(t)\|_{H^*(\Omega)}\le\|u(0)-v(0)\|_{H^*(\Omega)},
\end{equation}
which turns out to be also compact in $H^*(\Omega)$.
\end{thm}
  Note that the theorem applies to solutions with any sign.   The proof is based on showing that the fractional Laplacian operator can be suitably defined as a maximal monotone operator on the Hilbert space $H^*$. We follow Brezis's approach \cite{Brezis1} in dealing with  the standard Laplacian case where the operator is characterized as the sub-differential of a convex functional. The extention of this approach to the fractional Laplacian case is done in Proposition \ref{brezis}, after the needed functional analysis has been prepared in Section \ref{sec.exist}. The method produces not only existence and uniqueness of a semigroup of solutions, but also a number of important estimates, typical of evolution processes governed by maximal mononotone operators, cf. the classical monograph \cite{BrezisBk}. In particular, formula \eqref{strong} shows that solutions are actually $H^*$-strong in the sense of semigroup theory.  It is then easy to show that such solutions are indeed bounded weak energy solutions in the sense of \cite{DPQRV2}.  If moreover they are nonnegative, they are indeed strong $\LL^1$ solutions, in the sense that $u_t\in \LL^\infty((\tau,\infty): \LL^1(\Omega))$ for all $\tau>0$\,, see the remarks at the end of Section \ref{sec-Brezis}. This last property is convenient in the proof of the asymptotic behaviour.

\subsection{Asymptotic behaviour}

In Sections \ref{sec.asymp}, \ref{sec.fg}, and \ref{sec.entr} we discuss the second topic of the paper, namely large-time properties for solutions of \eqref{FPME.Problem.PHI} with $\varphi(u)=u^m$\,, with $m>1$. To this end it will be useful to introduce the following rescaled problem.
\begin{equation}\label{log.rescaling}
v(t,x)=(1+\tau)^{\frac{1}{m-1}}u(\tau,x)\,,\qquad t=\log(1+\tau)
\end{equation}
which transforms problem \eqref{FPME.Problem.PHI} into
\begin{equation}\label{FPME.prob.log}
\begin{split}
\left\{\begin{array}{lll}
v_{t}+ \A (v^m)=\dfrac{v}{m-1} & ~ {\rm in}~(0,+\infty)\times\Omega,\\[3mm]
v(0,x)=u_0(x) & ~{\rm in}~ \Omega, \\[3mm]
v(t,x)=0 & ~{\rm for}~  t >0 ~{\rm and}~ x\in\Gamma\,,
\end{array}\right.
\end{split}
\end{equation}
where the last line is a formal way of expression that we are considering fractional Laplacian operators with zero Dirichlet boundary conditions. Our first result is the existence of special separate-variable solution that we call the Friendly Giant following the denomination introduced by Dahlberg and Kenig for the standard
porous medium equation \cite{DahlKenig}. We also show its role as attractor of the evolution.

\begin{thm}\label{convLinf}
There exists a unique self-similar solution of the Dirichlet Problem \eqref{FPME.Problem.PHI} of the form
\begin{equation}\label{FG.intro}
U(\tau,x)=\frac{S(x)}{\tau^{\frac{1}{m-1}}}\,,
\end{equation}
for some bounded function $S: \Omega \to \RR$. Let $u$ be any nonnegative $H^*$-solution to the Dirichlet Problem \eqref{FPME.Problem.PHI}\,, then we have (unless $u\equiv 0$)
\begin{equation}\label{conv1}
\lim_{\tau\to \infty}\tau^{\frac{1}{m-1}}\left\|u(\tau,\cdot)-U(\tau,\cdot)\right\|_{\LL^\infty(\Omega)}=\lim_{t\to\infty}\left\|v(t,\cdot)-S\right\|_{\LL^\infty(\Omega)}=0\,.
\end{equation}
where $v$ is the solution of the rescaled flow \eqref{FPME.prob.log}\,.
\end{thm}
\noindent The previous theorem admits the following corollary.
\begin{cor}\label{thStat}
Let $m>1$.  For any given $c>0$  there exists a unique weak solution  to the elliptic equation
\begin{equation}\label{FPME.stat}
\A (S^m)=c\,{S} \qquad {\rm in}~\Omega,
\end{equation}
taking zero Dirichlet data in the sense defined above.
\end{cor}

\noindent Notice that the previous theorem is obtained in the present paper through a parabolic technique.

\noindent Once we have proven convergence of nonnegative solutions to the unique stationary state $S$, we look for sharp rates of convergence.

\begin{thm}\label{PME.Rates} Let $m>1$, let $v$ be  the rescaled solution as in \eqref{log.rescaling}, that converges to its unique
stationary state $S$. Then, we have
\begin{equation}
\|v(t,\cdot)-S(\cdot)\|_{\LL^1(\Omega)} \le\,K\,\ee^{-t}
\end{equation}
for all large $t\ge \overline{t}\gg 1$. It means that  as $\tau\to \infty$
\begin{equation}\label{conv1b}
\tau^{\frac{1}{m-1}}\left\|u(\tau,\cdot)-U(\tau,\cdot)\right\|_{\LL^\infty(\Omega)}=O(1/\tau)\,.
\end{equation}
\end{thm}
\noindent  The time rate of the error is sharp as can be seen by taking as $u$ a time-displaced version of $U$, i.e., $u(t,x)=U(t+t_1,x)$, so that
 $$
 u(t,x)-U(t,x)=-S(x)t^{-1/(m-1)}(1-(1+(t_1/t))^{-1/(m-1)}= S(x)(1-O(1/t)).
 $$
 In view of this observation we  can do better with respect to space.  Indeed, we can prove quantitative sharp rates of convergence for the relative error $w=u/U-1=v/S-1$, as follows.

\begin{thm}\label{thm.PME.Rates.Giant}
Let $u$ be any nonnegative $H^*$-solution to the Dirichlet Problem \eqref{FPME.Problem.PHI}\,, then we have (unless $u\equiv 0$) that there exist $t_0>0$ of the form
\[
t_0=\overline{k}\left[\frac{\int_{\Omega}\Phi_1\dx}{\int_{\Omega}u_0\Phi_1\dx }\right]^{m-1}
\]
such that for all $t\ge t_0$ we have
\begin{equation}\label{conv.rates.rel.err}
\left\|\frac{u(t,\cdot)}{U(t,\cdot)}-1 \right\|_{\LL^\infty(\Omega)} =\left\|\frac{v(t,\cdot)}{S(\cdot)}-1 \right\|_{\LL^\infty(\Omega)}\le \frac{2}{m-1}\,\frac{t_0}{t_0+t}\,.
\end{equation}
We remark that the constant $\overline{k}>0$ only depends on $m,d,s,$ and $|\Omega|$ and has explicit expressions given in the proof.
\end{thm}

\noindent\textbf{Remarks. }(i) The rate of convergence in these results is sharp, it can be checked in the same way as we have just done. See also the beginning of Section \ref{sec.fg}.

\noindent (ii) The convergence in the Laplacian case $s=1$ is a well-known result  by Aronson and Peletier \cite{Ar-Pe}.   Though we will follow their outline of proof, there are substantial difficulties due to the analysis of the boundary behaviour of solutions of the equations involving fractional operators.

\noindent (iii) We provide two different proofs of the convergence rates, namely Theorem \ref{PME.Rates} and Theorem \ref{thm.PME.Rates.Giant}. On one hand, Theorem \ref{thm.PME.Rates.Giant} gives a stronger convergence result and is somehow more explicit, but it heavily depends on the sharp estimates of \cite{BV2013}, cf. Theorem \ref{thm.GHP.PME}. On the other hand, Theorem \ref{PME.Rates} is obtained with a new entropy method that is based on the Caffarelli-Silvestre extension \cite{CaSi}, cf. also \cite{DPQRV2}, and it does not require strong estimates to get the sharp rate in $\LL^1$\,; in order to extend the sharp rates from $\LL^1$ to $\LL^\infty$ norm, we would need $C^\alpha$ regularity.

\subsection{Sobolev spaces and interpolation. }
In Appendix A.1 we review in some detail for the readers' convenience the theory of Fractional Sobolev Spaces which is needed to characterize our fractional Laplacian operators. The second appendix describes the interpolation method known as discrete J-Method, which we use to identify the spaces $H$, defined trough eigen-elements of $\A$\,, with the fractional Sobolev spaces $H^s$, usually constructed via interpolation.

%%%%%%%%%%%%%%%%%%%%%%%%%%%%%%%%%%%%%%%%%%%%%%%%%%%%%%%%%%%%%%%%%%%%%%%%%%%%%%%%%%%%
\section{Existence and uniqueness of $H^*$-solutions}\label{sec.exist}

In order to construct solutions of Problem \eqref{FPME.Problem.PHI} we are going to use the theory of semigroups in Hilbert spaces generated by the subdifferentials of convex functions, following the ideas of Brezis \cite{Brezis1}, where the standard Laplacian case is treated,  see also the general theory  of \cite{BrezisBk} and the account of \cite{VazBook} for the PME\,. This part of the theory can be done for a quite general class of linear operators $\A$ and nonlinearities $\varphi$  as we will comment later.

%%%%%%%%%%%%%%%%%%%%%%%%%%%%%%%%%%%%%
\subsection{Functional setting}

 As a first step, we clarify the different functional spaces involved in the theory. More precisely, we want to identify the functional space $H(\Omega)$, corresponding to the each of our operators $\A$ according to formula \eqref{defH}, in terms of the more standard Sobolev spaces $H^s(\Omega)=W^{s,2}(\Omega)$. This is crucial for our purposes because we want to use the abstract monotone operators setting, and to do that we need Sobolev-type inequalities and compact imbedding of $H$ into $L^p$ spaces. Since this relies on known theories, we just sketch hereafter the main points, more complete details can be found in Appendix 1.

%%%%%%%%%%%%%%%%%%%%%%%
\subsubsection{Spaces $H$ and $H^*$ associated to a self-adjoint operator}\label{sect.H.spaces}

Let $ \A : \dom( \A )\subset\LL^2(\Omega)\to \LL^2(\Omega)$ be a  positive self-adjoint operator, where $\Omega$ is a bounded domain with smooth boundary as above. Assume that $ \A $ has a discrete spectrum and a $\LL^2(\Omega)$ ortonormal basis of eigenfunctions. We denote by $\lambda_k$  its eigenvalues written in increasing order and repeated according to their multiplicity, and by  $\{\Phi_k\}$ the corresponding $\LL^2(\Omega)$-normalized eigenfunctions; they form an ortonormal basis for $\LL^2(\Omega)$\,.    In this generality we are also assuming that $0$ is not an eigenvalue, which is true for Dirichlet problems.   We can always associate to $ \A $ the bilinear form
\[
B(f,g):=\int_\Omega f\,  \A  g\dx = \int_\Omega g\,  \A  f\dx= \int_\Omega \AM g\,  \AM  f\dx\,,
\]
which turns out to be a Dirichlet form whose completed domain we call $H\times H \subset \LL^2(\Omega)\times\LL^2(\Omega)$\,. Define the norm:
\begin{equation}\label{def.H.norm}
\|f\|_H=\left(\sum_{k=1}^{\infty}\lambda_k\hat{f}_k^2\right)^{\frac{1}{2}}<+\infty\qquad\mbox{with}\qquad \hat{f}_k=\int_\Omega f(x)\Phi_k(x)\dx\,.
\end{equation}
It is standard to show that the closure of the domain of the Dirichlet form of $ \A $ is given by
\begin{equation}\label{def.H.space}
H=H(\Omega):=\left\{f\in \LL^2(\Omega)\;|\; \sum_{k=1}^{\infty}\lambda_k\hat{f}_k^2<+\infty\right\}
\end{equation}
 in other words, $H$ is the domain of $\AM$.   The above function space  is a Hilbert space with the inner product given by the Dirichlet form
\[
\langle f, g \rangle_H = \sum_{k=1}^{\infty}\lambda_k\hat{f}_k\,\hat{g}_k=B(f,g)
\]
(We sketch a proof of the above facts as a consequence of integration by parts in  Appendix \ref{Proofs.Sect.211}).

\noindent We will also consider the dual space $H^*$ endowed with its dual norm,
\begin{equation}\label{def.H*.norm.1}
\|F\|_{H^*}=\sup_{\begin{subarray}{c}  g\in H \\ \|g\|_{H}\le 1\end{subarray}}\langle F, g \rangle_{H^*,H}
=\sup_{\begin{subarray}{c}  g\in H \\ \|g\|_{H}\le 1\end{subarray}}\sum_{k=1}^{\infty}\hat{F}_k\,\hat{g}_k\,,
\end{equation}
where $\langle \cdot\,,\, \cdot \rangle_{H^*,H}$ is the duality mapping. We have
\begin{equation}\label{def.H*.norm.2}
\|F\|_{H^*}=\left(\sum_{k=1}^{\infty}\lambda_k^{-1}\hat{F}_k^2\right)^{\frac{1}{2}}\,.
\end{equation}
It is clear at this point that $ \A $ gives the canonical isomorphism between $H$ and $H^*$\,.

\medskip

\noindent\textbf{Integration by parts in $H$. }The Spectral Theorem allows to write $ \A $ as
\[
 \A  f(x)=\sum_{k=1}^{\infty}\lambda_k \hat{f}_k\,\Phi_k(x)\qquad\mbox{for any $f\in H$\,,}
\]
Therefore, the integration by parts formula follows:
\begin{equation}\label{Int.by.Parts}
\int_\Omega f  \A  g\dx=\int_\Omega \A^{1/2} f \,  \A^{1/2} g\dx=\int_\Omega g  \A  f\dx
\qquad\mbox{for any $f,g\in H$.}
\end{equation}

\medskip

\subsubsection{The inverse operator $\AI$ and the Green function}

Define the inverse operator: $\AI: H^*(\Omega)\to H(\Omega)$:
\begin{equation}\label{Green.Operator.Definition}
\AI f(x_0)=\int_\Omega \K(x,x_0)f(x)\dx
\end{equation}
which is the inverse of the canonical isomorphism $\A:H(\Omega)\to H^*(\Omega)$:
let $f\in H^*(\Omega)$ and $F$ be the solution to the Dirichlet problem
\begin{equation}\label{Green.Operator.Equation}
\left\{\begin{array}{lll}
\A(F)=f &  ~ {\rm in}~ (0,+\infty)\times \Omega\\
F(x)=0 & ~\mbox{on $(0,\infty)\times\Gamma$}
\end{array}
\right.
\end{equation}
When $\A$ is either the SFL of the RFL, the kernel of the inverse can be expressed in terms of the Green function $\K(x,x_0):=G_\Omega(x,x_0)$ and satisfies the following estimates, if the domain $\Omega$ is smooth enough:
\begin{equation}\label{Gree.est.1}
c_{0,\Omega}\Phi_1(x)\Phi_1(y)\le \K(x,y)\le \frac{c_{1,\Omega}}{|x-x_0|^{N-2s}}
\left(\frac{\Phi_1(x)}{|x-x_0|^\gamma}\wedge 1\right)
\left(\frac{\Phi_1(x_0)}{|x-x_0|^\gamma}\wedge 1\right)
\end{equation}
with $\gamma=1$ for the SFL and $\gamma=s$ for the RFL, where $\Phi_1$ is the first eigenfunction of $\A$, which satisfies the following estimates
\begin{equation}\label{Phi1.est.1}
k_{0,\Omega}\, \big(\dist(x, \partial\Omega)^\gamma\wedge 1\big) \le \Phi_1(x)
\le k_{1,\Omega}\, \big(\dist(x, \partial\Omega)^\gamma\wedge 1\big)\,.
\end{equation}
 The above estimates \eqref{Gree.est.1} can be found in \cite{Jak, Kul}, for the RFL, while in the case of the SFL they can be easily deduced by the celebrated heat kernel estimates of Davies and Simon \cite{Davies1, DS} (for the case $s=1$), see also \cite{BV2013} for more details about the application of the results of \cite{Davies1, DS} to the SFL.

\noindent Finally, we recall that $\AI$ is an isomorphism of the two Hilbert spaces $H^*(\Omega)$ and $H(\Omega)$ as explained in the previous paragraph.

\noindent\textbf{Integration by parts in $H^*$. }The Spectral Theorem allows to write $\AI$ as the series
\[
\AI f(x)=\sum_{k=1}^{\infty}\lambda_k^{-1} \hat{f}_k\,\Phi_k(x)\qquad\mbox{for any $f\in H^*$\,.}
\]
therefore the integration by parts formula
\begin{equation}\label{Int.by.Parts.H*}
\int_\Omega f \AI g\dx=\int_\Omega \AIM f \,  \AIM g\dx=\int_\Omega g \AI f\dx
\qquad\mbox{for any $f,g\in H^*$\,,}
\end{equation}
follows by the orthogonality property of the eigenfunctions $\Phi_k$, exactly as for the operator $\A$\,.
As a consequence of this formula, we can express the scalar product and the norm of $H^*$ in the following equivalent ways, for any $f,g\in H^*$\:
\begin{equation}\label{product.H*}
\langle f, g \rangle_{H^*} = \sum_{k=1}^{\infty}\lambda_k^{-1}\hat{f}_k\,\hat{g}_k
=\int_\Omega \AIM f \,  \AIM g\dx
=\int_\Omega f \AI g\dx=\int_\Omega g \AI f\dx\,,
\end{equation}
hence
\begin{equation}\label{norm.H*}
\|f\|_{H^*} = \sum_{k=1}^{\infty}\lambda_k^{-1}\hat{f}_k^2
=\left\|\AIM f\right\|_{\LL^2(\Omega)}^2
=\int_\Omega f \AI f\dx\,.
\end{equation}

%%%%%%%%%%%%%%%%%%%%%%%%%%%%%%%%%%%%%%%%%%%%%%%%%%%%%%%%%%%%%%%%%%%%
\subsubsection{Relation between  the $H$ space and  $H^s(\Omega)=W^{2,s}(\Omega)$ spaces }\label{identif.H.Hs}

Even if we can deal with a quite wide class of linear operators, we concentrate our attention from here on, to the cases where $\A$ is either the Spectral fractional Laplacian (SFL) or the Restricted Fractional Laplacian (RFL). We are going to identify the corresponding space $H$ with the fractional Sobolev spaces; namely, we are going to prove that both for the SFL and the RFL we have
\begin{equation}\label{H-Hs}
H(\Omega)=\left\{\begin{array}{lll}
H^s_0(\Omega)\,, &\qquad\mbox{if }\frac{1}{2}<s\le 1\,,\\[2mm]
H^{1/2}_{00}(\Omega)\,,&\qquad\mbox{if }s=\frac{1}{2}\,,\\[2mm]
H^s(\Omega)\,, &\qquad\mbox{if }0<s< \frac{1}{2}\,.
\end{array}\right.
\end{equation}
 The precise definition of these spaces will be given in Section \ref{sec.7.3}\,. We recall that it is possible to obtain a unified characterization of them\,, namely \eqref{H-Hs} can be written as $H(\Omega)= \dot H^s( \overline{\Omega}) = \{\,u \in H^s(\RR^d)\;|\;\supp(u) \subset \overline{\Omega}\,\}$\,, cf. \cite{McLean, Triebel}.

We give next the proof of the identification \eqref{H-Hs} for both operators under consideration. It relies on interpolation theory, more precisely on the Discrete version of J-Method of Theorem \ref{Thm.Jdiscr.Appendix}, proved in Appendix 2.

\medskip

\noindent\textbf{Case I. The spectral fractional Laplacian}

\noindent Let  $\A = (-\Delta_\Omega)^s$ be the spectral fractional Laplacian, as in \eqref{sLapl.Omega.Spectral}. We have already denoted the eigenfunctions and eigenvalues  by $(\lambda_k^s,\phi_k)$, where $\lambda_k$ are the eigen-elements of the Dirichlet Laplacian (i.e., the classical case $s=1$), so that
\begin{equation}\label{uk.spectral}
u(x)=\sum_{k\ge 1}u_k(x)=\sum_{k\ge 1}\hat{u}_k\phi_k(x)\quad\mbox{with}\quad \hat{g}_j=\int_\Omega g(x)\phi_j(x)\dx\,,\quad\mbox{and}\quad \|\phi_j\|_{\LL^2(\Omega)}=1\,.
\end{equation}
Recall that $\|u_k\|_{\LL^2(\Omega)}=\hat{u}_k \|\phi_k\|_{\LL^2(\Omega)}=\hat{u}_k$\,, so that
\[
\|u_k\|_{H^1_0(\Omega)}^2=\hat{u}_k^2\int_\Omega |\nabla \phi_k|^2\dx=\hat{u}_k^2\int_\Omega \phi_k\Delta\phi_k\dx=\lambda_k\hat{u}_k^2\,.
\]
Let now $0<s<1$ with $s\ne 1/2$, and recall that $H^s_0$ can be defined by interpolation as in \eqref{def.Hs0.Omega.2}, namely
\begin{equation*}
H^s_0(\Omega):= \big[H^{s_1}_0(\Omega)\,,\,H^{s_2}_0(\Omega)\big]_\theta\,,\qquad\mbox{for any $0<s_1\,,\, s_2\le 1$ such that $(1-\theta)s_1+\theta s_2\neq\frac{1}{2}$\,.}
\end{equation*}
We recall the definition \eqref{def.H.space} of $H$, namely
\begin{equation*}
H:=\left\{f\in \LL^2(\Omega)\;|\; \sum_{k=1}^{\infty}\lambda_k^s\hat{f}_k^2<+\infty\right\}
\end{equation*}Therefore, in order to apply the discrete version of the J-Method for interpolation, see Appendix, for $u\in H$  we  define the sequence
\begin{equation}
U_k=\lambda_k^{-\theta} J(\lambda_k,u_k)
=\lambda_k^{-\theta}\max\{\|u_k\|_{H^1_0(\Omega)}\,,\,\lambda_k\|u_k\|_{\LL^2(\Omega)}\}=\lambda_k^{1-\theta}\hat{u}_k
\end{equation}
and it is clear that $U_k\in \ell^2(\NN)$ if and only if $\theta=1-s$ and $u\in H$\,,
and that $\|u\|_H=\|U_k\|_{\ell^2(\NN)}$\,. The discrete version of the J-Method, namely Theorem \ref{Thm.Jdiscr.Appendix}, allows to identify $H$ as the interpolation space:
\[
H=\big[H^1_0(\Omega)\,,\,\LL^2(\Omega)\big]_{1-s}=H^s_0(\Omega)\,.
\]
and the norms on $H^s_0(\Omega)$ and $H$ are equivalent. Finally, when $s=1/2$\,, the above discussion can be repeated, but we identify $H$ as the interpolation space:
\[
H=\big[H^1_0(\Omega)\,,\,\LL^2(\Omega)\big]_{1/2}=H^{\frac{1}{2}}_{00}(\Omega)\,.
\]
and again the norms on $H^{\frac{1}{2}}_{00}(\Omega)$ and $H$ are equivalent.

\noindent\textbf{Remark. }In order to apply Theorem \ref{Thm.Jdiscr.Appendix}, one has to check a growth condition $0<\lambda_{k+1}/\lambda_k<\Lambda_0<+\infty$, for some $\Lambda_0>0$. In the case of the SFL, this property is satisfied since we know that $\lambda_k\asymp k^{2/d}$, see e.g. \cite{Davies2}\,.

\medskip

\noindent\textbf{Case II: The restricted fractional Laplacian}

\noindent Let $\A = (-\Delta_{|\Omega})^s$ be the RFL, corresponding to the expression in terms of hypersingular kernel \eqref{sLapl.Rd.Kernel} restricted to functions supported in $\Omega$, with external conditions \eqref{FPME.Dirichlet.conditions.Restricted}. The eigen-elements of $\A$ are $(\lambda_{s,k},\phi_{s,k})$ and it is known that $\lambda_{s,k}\le \lambda_k^s$,    where $\lambda_k^s$ are the eigenvalues of the SFL, see \cite{ChSo} for a proof. Moreover, the eigenfunctions are different from the eigenfunctions of the SFL, as it has already been discussed in the Introduction. Therefore, for $u\in \LL^2(\Omega)$:
\begin{equation}\label{uk.restricted}
u(x)=\sum_{k\ge 1}u_k(x)=\sum_{k\ge 1}\hat{u}_k\phi_{s,k}(x)\quad\mbox{with}\quad \hat{g}_j=\int_\Omega g(x)\phi_{s,j}(x)\dx\,,\quad\mbox{and}\quad \|\phi_{j,s}\|_{\LL^2(\Omega)}=1\,.
\end{equation}
Recall that $ \|u_k\|_{\LL^2(\Omega)}=\hat{u}_k \|\phi_{s,k}\|_{\LL^2(\Omega)}=\hat{u}_k$\,.
The proof continues in quite different way from SFL, since our eigenfunctions are not related to the standard Laplacian anymore.

We observe that  $\phi_{s,k}\in H^{2s}_0(\Omega)$. Indeed, recall that $E_0(\phi_{s,k})$ is the extension by zero outside $\Omega$ of $\phi_{s,k}$.  Therefore,
\[
\begin{split}
\|u_k\|_{H^{2s}_0(\Omega)}^2
&=\hat{u}_k^2\|\phi_{s,k}\|_{H^{2s}_0(\Omega)}^2
=\hat{u}_k^2 \|E_0(\phi_{s,k})\|_{H^{2s}(\RR^d)}
=\hat{u}_k^2 \|(-\Delta_{\RR^d})^s\,E_0(\phi_{s,k})\|_{\LL^2(\RR^d)}^2\\
&=\lambda_{s,k}^2\hat{u}_k^2 \|E_0(\phi_{s,k})\|_{\LL^2(\RR^d)}^2
=\lambda_{s,k}^2\hat{u}_k^2 \|\phi_{s,k}\|_{\LL^2(\Omega)}^2=\lambda_{s,k}^2\hat{u}_k^2
\end{split}
\]
where we have taken as a norm on $H^{2s}_0$ the equivalent norm defined through zero extension, as in \eqref{Hs.Omega.norm.2}\,.

Let $s\ne 1/2$, and recall that $H^s_0$ can be defined by interpolation as in \eqref{def.Hs0.Omega.2}, namely
\begin{equation*}
H^s_0(\Omega):= \big[H^{s_1}_0(\Omega)\,,\,H^{s_2}_0(\Omega)\big]_\theta\,,\;\mbox{for any $0\le s_1\,,\, s_2 $ such that $s=(1-\theta)s_1+\theta s_2\neq \mbox{integer}+\frac{1}{2}$\,.}
\end{equation*}
Next we would like to apply the discrete version of the J-Method, namely Theorem \ref{Thm.Jdiscr.Appendix} and to this end we define the sequence
\begin{equation}
U_k=\lambda_k^{-\theta} J(\lambda_k,u_k)
=\lambda_k^{-\theta}\max\{\|u_k\|_{H^{2s}_0(\Omega)}\,,\,\lambda_k\|u_k\|_{\LL^2(\Omega)}\}=\lambda_k^{1-\theta}\hat{u}_k
\end{equation}
and it is clear that $U_k\in \ell^2(\NN)$ if and only if $\theta=1/2$ since $u\in H$\,, where $H$ is the domain of the associated Dirichlet form, as discussed in Section \ref{sect.H.spaces}, namely
\begin{equation*}
H:=\left\{f\in \LL^2(\Omega)\;|\; \sum_{k=1}^{\infty}\lambda_k\hat{f}_k^2<+\infty\right\}
\end{equation*}
and it is clear that $\|u\|_H=\|U_k\|_{\ell^2(\NN)}$\,. To apply the J-method, we need bounds on the eigenvalues for $\A$, the restricted Laplacian. For that, we invoke the following result in \cite{BG}: the eigenvalues $\lambda_{s,k}$ of the RFL behave like
$$\lambda_{s,k}=\frac{d+2s}{d}C_{d,s}|\Omega|^{-2s/d}k^{2s/d}(1+O(1))\,,\;\mbox{as }k \to \infty. $$
Then one has that
$$0<\lambda_{s,k+1}/\lambda_{s,k}<\Lambda_0<+\infty,$$
hence we can apply the discrete version of the J-Method, namely Theorem \ref{Thm.Jdiscr.Appendix}, which implies that
\[
H=\big[H^{2s}_0(\Omega)\,,\,\LL^2(\Omega)\big]_{1/2}=H^s_0(\Omega)\,.
\]
and the norms on $H^s_0(\Omega)$ and $H$ are equivalent.

\noindent Finally, when $s=1/2$\,, the above discussion can be repeated verbatim, but we identify $H$ as the interpolation space:
\begin{equation}\label{interp.0}
H=\big[H^{1}_0(\Omega)\,,\,\LL^2(\Omega)\big]_{1/2}=H^{\frac{1}{2}}_{00}(\Omega)\,.
\end{equation}
and the norms on $H^{\frac{1}{2}}_{00}(\Omega)$ and $H$ are equivalent.  Notice that only when $s=1/2$ the interpolation spaces and exponents in formula \eqref{interp.0} are exactly the same both for the SFL and the RFL.

 \noindent\textbf{A direct proof of the previous result for $s \neq \frac12$}
We estimate
$$
\int_\Omega f(x) \A f(x)\,dx=\int_\Omega f(x)\,dx \int_{\RR^d}\frac{f(x)-f(y)}{|x-y|^{d+2s}}\,dy.
$$
Hence
$$
\int_\Omega f(x) \A f(x)\,dx=\int_{\RR^d} \int_{\RR^d} E_0(f)(x)\frac{E_0(f)(x)-E_0(f)(y)}{|x-y|^{d+2s}}\,dx \,dy.
$$
By symmetry, this leads to
$$
\|f\|^2_H=\int_\Omega f(x) \A f(x)\,dx=C\int_{\RR^d} \int_{\RR^d}\frac{|E_0(f)(x)-E_0(f)(y)|^2}{|x-y|^{d+2s}}\,dx \,dy=C\|E_0(f)\|^2_{W^{s,2}(\RR^d)}.
$$
The previous computation holds for any $C^\infty_0(\Omega)$ function. By approximation, it also holds for functions in $\LL^2(\Omega)$. Therefore, the only point is to identify the semi-norm $\|E_0(f)\|_{W^{s,2}(\RR^d)}$ with the equivalent $H^s(\Omega)$ norm, to obtain formula \eqref{H-Hs}. This can be done only for $s\neq 1/2$, since the extension by zero outside $\Omega$ is not continuous from $H^{1/2}(\Omega)\to H^{1/2}(\RR^d)$, see Section \ref{sect.prop.Hs} for further details; roughly speaking we can say that $H^{1/2}$ is not a restriction space.
When $s\ne 1/2$, we can identify $\|E_0(f)\|_{W^{s,2}(\RR^d)}$ with the equivalent $H^s(\RR^d)$ norm of $E_0(f)$ (see also Theorem \ref{thm.Hs-Ws}), which is equivalent to the norm of $H^s(\Omega)$, see \eqref{Hs.Omega.norm.2}.

\noindent\textbf{Summary of the previous discussion: }Let $\A$ be one of the two operators under consideration, i.e. the SFL or the RFL. Then $\A:H\to H^*$ is an isomorphism between the space $H$ (closure of the domain of the Dirichlet form associated to $\A$) and its dual $H^*$. Moreover, $H$ can be characterized as in \eqref{H-Hs}\,.

\noindent\textbf{Case III: Other examples}\label{more.examples}

Though the existence theory we develop in the next section has been devised with these two operators in mind, it is in fact quite general. Indeed, we can handle any operator having a discrete set $\left \{\lambda_k, \Phi_k \right \}$ of eigen-elements, $\{\Phi_k\}$ being an orthonormal basis in $\LL^2(\Omega)$\,.  The sole further assumption that is required to identify the space $H$ with the fractional Sobolev spaces $H^{s}_0$ is that the eigenvalues do not grow more than exponentially. Let us give hints of interesting cases where the results of the paper apply.

\noindent $\bullet$\textit{Spectral type: } In particular, one could replace the Laplacian by operators $ \mathcal L$ being self-adjoint in $L^2(\Omega)$ and positive given by powers of differential operators of the type
$$
{\mathcal A}= \mbox{div}(A\nabla)
$$
with homogeneous Dirichet conditions and $A$ uniformly elliptic, sufficiently smooth, i.e., an elliptic operator with coefficients.
Then we can follow the spectral type above and consider the fractional powers of this operator.\\
\sl Existence: \rm the eigenvalues estimates $\lambda_k\asymp k^{2s/d}$ can be obtained as in \cite{Davies2} Theorem 6.3.1\. \\
\sl Asymptotics: \rm from the strong Heat-kernel estimates of  \cite{Davies1,DS}, we can derive Green function estimates. In \cite{BV2013} we show how such heat kernel bounds give the estimates of Theorem \ref{thm.GHP.PME}.

\noindent $\bullet$\textit{Restricted type: }Integral operators with a kernel which has a singularity of the type $|x-y|^{d+2s}$ can be treated, under some further assumptions. Clearly, one condition is that the eigenvalues do not grow more than exponentially.  Another condition would be:
\[
\int_{\RR^d} f(x)\A f(x) \dx \asymp \int_{\RR^d}\int_{\RR^d}\frac{|f(x)-f(y)|^2}{|x-y|^{d+2s}}\dx\dy\qquad\mbox{for all }f\in C_c^\infty(\Omega)\,.
\]
The existence and uniqueness theory works at least under one of the above additional conditions.

\medskip

\subsection{Existence and uniqueness results. Proof of Theorem \ref{cor3.2}}\label{sec-Brezis}

The proof of Theorem  \ref{cor3.2} follows essentially by applying the techniques of monotone operators in Hilbert spaces. More precisely our Theorem  \ref{cor3.2} is the $H^*$ version of Brezis' result in $H^{-1}(\Omega)$, namely Corollary 31 of his paper \cite{Brezis1}. The proof can be easily adapted to our operators once the proper functional analysis is in place. Next, we will review the main steps  for the reader's convenience. Let us begin by setting up some notations.

Let $j$  be a convex, lower semi-continuous function $j:\RR\to \RR$ and such that $j(r)/|r|\to \infty$ as $|r|\to\infty$. We let $\varphi=\partial j$ be the sub-differential of $j$.  For $u\in H^*(\Omega)$ we define
\[
\Psi(u)=\int_\Omega j(u)\,\dx
\]
whenever $u\in\LL^1(\Omega)$ and $j(u)\in\LL^1(\Omega)$, and define $\Psi(u)=+\infty$ otherwise. The example we have in mind is  $\varphi(u)=|u|^{m-1}u$ and $j(u)=|u|^{m+1}/(m+1)$\,, so that $\Psi(u)=\|u\|_{\LL^{m+1}(\Omega)}/(m+1)$\,. Notice that here we can consider a wider range of powers, namely any $m>0$.

 \noindent Staying   in the more general case of $j$ and $\varphi$, we need to prove the following Proposition.

\begin{prop}\label{brezis}
 $\Psi$ is convex and lower semi-continuous function in $H^*(\Omega)$, so that its sub-differential $\partial\Psi$ is a maximal monotone operator in $H^*(\Omega)$. Moreover, this sub-differential $\partial\Psi$ can be characterized as follows:
\begin{equation}
f\in\partial\Psi(u) \qquad\mbox{if and only if}\qquad \AI f(x)\in\varphi(u(x))\qquad\mbox{ a.e. in }\Omega.
\end{equation}

\end{prop}
\noindent {\sl Proof.~}  Following Brezis' arguments, the proof of Theorem \ref{cor3.2}  will be split in several steps. It is worth recalling here that an essential point in the proof is that we can identify the space $H$ with the usual Sobolev spaces, see formula \eqref{H-Hs}\,, hence the usual Sobolev imbeddings and inequalities of $H^s_0(\Omega)=W^{2,s}_0(\Omega)$ in $\LL^p$ spaces hold for $H$ as well\,; we refer to the Appendix A1 for further details on this issues.

\noindent$\bullet~$\textsc{Step 1. }We first need to prove that the functional  $u\mapsto \Psi(u)$ is l.s.c. on $\LL^1(\Omega)$. By Fatou's Lemma it is sufficient to see that this functional is convex l.s.c. on $\LL^1(\Omega)$\,. The proof of this statement is exactly the same in the proof of Thm. 17 of \cite{Brezis1}, pg. 123-124, just by replacing $H^{-1}$ with $H^*$.

We define the operator $A$ on $H^*$ to be
\[
A u=\left\{\A w\;|\; w\in H\mbox{ and }w(x)\in \varphi(u(x))\mbox{ a.e. on }\Omega \right\}
\]
with $u\in \dom(A)$ if and only if there is some $w\in H$ such that $w\in\varphi(u(x))\mbox{ a.e. on }\Omega$.
We want to prove that $A\subset\partial\Psi$ (Step 2) and then that $A$ is maximal monotone (Step 3). Before doing that, we need an approximation result that we state and prove in Step 2

\noindent$\bullet~$\textsc{Step 2. } We prove first that $A\subset\partial\Psi$. This is formulated as follows: Let $f\in Au$, i.e. $u \in H^*\cap \LL^1(\Omega)$\,, $f=\A w$ with $w\in H$\,, $w(x)\in \varphi(u(x))$ a.e. on $\Omega$\,. We need to prove that under these assumptions $\Psi(u)$ is bounded and for every  $v \in H^*\cap \LL^1(\Omega)$ be such that $j(v)\in \LL^1(\Omega)$ we have
$$
\Psi(v)-\Psi(u)\ge \langle f,v-u\rangle_{H^*\times H^*}=\langle w,v-u\rangle_{H\times H^*}\,,
$$
where we have used the identification of $H$ with $H^*$ via the isomorphism $\cal L$. Therefore, we need
to prove that
\begin{equation}\label{ineq.A}
 \int_\Omega (j(u)-j(v))\dx \le \langle w,u-v\rangle_{H\times H^*}\,.
\end{equation}
Now, we have that $j(u)-j(v)\le w (u-v)$  a.e. on $\Omega$\,. So the remaining question depends on the correct interpretation of the right-hand side and its comparison with the Lebesgue integral.
We state separately the needed result as a \\
{\bf Claim:} {\it Let $F\in H^*\cap \LL^1(\Omega)$ and let $w\in H.$ Let $g\in \LL^1(\Omega)$ and let $h$ be measurable with
\begin{equation}\label{Bre.2.1}
g\le h\le F\, w\qquad\mbox{a.e. on }\Omega\,.
\end{equation}
Then $h\in \LL^1(\Omega)$ and
\begin{equation}
\int_\Omega h(x)\dx\le \int_\Omega F(x)\,w(x)\dx\,.
\end{equation} }

\noindent To end of Step 2 we apply now the Claim  with $F=u-v$\,, $h=j(u)-j(v)$ and $g=-c_1|u|-c_2 -j(v)$\,, (with $j(r)\ge -c_1|r|-c_2$). We conclude that $j(u)\in \LL^1(\Omega)$ and
\eqref{ineq.A} holds. Hence, $f\in \partial\Psi(u)$\,.

\noindent {\sl Proof of the Claim.~}\rm Let
\[
w_n=\left\{
\begin{array}{lll}
n\,, &\qquad\mbox{if }w\ge n\,, \\
w\,, &\qquad\mbox{if }|w|\le n\,, \\
-n\,, &\qquad\mbox{if }w\le -n\,.
\end{array}
\right.
\]
Note that $|w_n/w|\le 1$ and $w_n/w\to 1$ a.\,e.in $\Omega$. On the other hand, we have the pointwise inequality: $|w_n(x)-w_n(y)| \leq |w(x)-w(y)|$ for every $x,y \in \Omega$. This
implies that the sequence $\{w_n\}$ is uniformly bounded in $H$, hence it weakly converges in $H$. Moreover, $w_n\to w$ a.e. in $\Omega$, then by dominated convergence we have that $w_n\to w$ strongly in $L^p$\,, for some $p\ge 1$ and also that $\|w-w_n\|_H\to 0$, i.\,e.,  $w_n \to w$ in $H$.

\noindent Let also $h_n=h\,w_n/w$ and let $g_n=g\,w_n/w$. Multiplying \eqref{Bre.2.1} by $w_n/w$
\[
g_n\le h_n\le F\, w_n\qquad\mbox{a.e. on }\Omega\,,
\]
and hence
\[
0\le h_n-g_n \le F\,w_n-g_n\qquad\mbox{a.e. on }\Omega\,.
\]
The sequence $h_n-g_n\to h-g$ a.e. on $\Omega$ as $n\to \infty$ and also
\[
\int_\Omega(h_n-g_n)\dx\le \int_\Omega F\,w_n \,\dx - \int_\Omega g_n\dx\,.
\]
Note that since $w_n\in L^2(\Omega)$ we have $\int_\Omega F\,w_n \,\dx=\langle F, w\rangle_{H^*\times H}$.
Since $w_n\to w$ in $H$ and $g_n\to g$ in $\LL^1(\Omega)$\,, we conclude by Fatou's lemma  that $h-g\in \LL^1(\Omega)$\,, and thus that $h\in  \LL^1(\Omega)$ with
\[
\int_\Omega(h-g)\dx\le \int_\Omega F\,w \,\dx - \int_\Omega g\dx\,.\mbox{\qed}
\]

\noindent$\bullet~$\textsc{Step 3. }We prove now that $A$ is maximal monotone. For a given $f\in H^*$ we have to find $u \in H^*\cap \LL^1(\Omega)$ and $w\in H$ such that $u+\A w= f$ and $w(x)\in \varphi(u(x))$ a.e. on $\Omega$\,. Let $\eta=\varphi^{-1}$ so that $\dom(\eta)=\RR$. Without loss of generality we can assume that $0\in\eta(0)$. Let $w_\mu\in H$ be the solution of the equation
\begin{equation}\label{Bre.4.1}
\eta(w_\mu)+\A w_\mu=f\,,
\end{equation}
which exists by standard results.  Multiplying \eqref{Bre.4.1} by $w_\mu$ and integrating over $\Omega$, we see that $w_\mu$ is bounded in $H$ as $\mu\to 0$. Thus we can find a sequence $\mu_n\to 0$ such that $w_{\mu_n}\to w$ weakly in $H$\,, $w_{\mu_n}\to w$ a.e. on $\Omega$\,, $(I+\mu_n\eta)^{-1}w_{\mu_n}\to w$ a.e. on $\Omega$.

\noindent The proof is concluded once we use the following result:

\it Let $\eta$ be a maximal monotone graph on $\RR\times\RR$ such that $\dom(\eta)=\RR$ and $0\in \eta(0)$\,. Let $f_n$ and $v_n$ be measurable functions on $\Omega$ such that $v_n \to v$ a.e. on $\Omega$\,, $f_n(x)\in \eta(v_n(x))$ a.e. on $\Omega$ and $f_n\,v_n\in \LL^1(\Omega)$ with $\int_\Omega f_n\, v_n\dx\le C$. Then there exists a subsequence $n_k\to +\infty$ such that $f_{n_k}\to f$ weakly in $\LL^1$ and $f(x)\in \eta(v(x))$ a.e. on $\Omega$\,.\rm

\noindent This is exactly Theorem 18 of \cite{Brezis1} pg. 126\,, which holds in general: the proof does not rely on properties of $H^1_0$, it only needs a Lemma (Lemma 3 of  \cite{Brezis1} pg. 126-127) which holds for maximal monotone operators on Hilbert spaces.   Corollary 31 of \cite{Brezis1} provides the result for $s=1$\,, but the same proof also holds for $0<s<1$\,, because it relies on the above results.

It remains to prove that the solution map $S_t: u_0\mapsto u(t)$ defines a semigroup of (non-strict) contractions in $H^{*}(\Omega)$, i.\,e., satisfying \eqref{contractivity.H*}. This holds in a wider generality and follows from remarks after Theorem 21 of \cite{Brezis1}, for the contractivity estimates see also Remark (i) below. The compactness in $H^*$ follows by the boundedness of $u(t)$ in the $H$ norm, together with compact Sobolev embeddings of $H$ into $H^*$. See also Theorems 6.17  and 6.18 of Chapter 6 of \cite{VazBook}.\qed

\noindent\textbf{Remarks. }We address the attention to the following two consequences of Theorem \ref{cor3.2}:\\
\noindent(i) \textit{Uniqueness. }Some interesting estimates hold for a.e. $t$ and in distributional sense:
\begin{equation}
\frac{\rd}{\dt}\|u_1-u_2\|^2_{H^*}= -2 \int  (u_1-u_2)(\varphi(u_1)-\varphi(u_2))\,dx,
\end{equation}
which immediately implies uniqueness and the contraction estimate \eqref{contractivity.H*} in $H^*$. We refrain from further details because similar computations have been done in the case $s=1$ in  Brezis' paper \cite{Brezis1}, see also Chapters 6, 10 and 11 of \cite{VazBook}.

\noindent(ii) \textit{Comparison.} Besides, the Comparison Principle holds in the sense that for two solutions $u,v$ with initial data $u_0\le v_0$ a.e. in $\Omega$, then $u\le v$ a.e. in $(0,\infty)\times \Omega$.   In particular, $u_0\ge 0$ implies $u(t)\ge 0$ for all $t>0$. All of this can be proved just as in the case $s=1$ see Chapter 3 of \cite{VazBook}.

\noindent(iii) \textit{$H^*$ solutions are energy solutions. }Since $t u_t\in H^*$ then $\varphi(u(t))\in H$ for almost all $t>0$ (in fact all), therefore $u$ is a weak energy solution in the sense of Definition 4.1 of \cite{DPQRV2}. Moreover, by comparison with weak energy solutions to the Cauchy problem in the whole space $\RR^d$\,, we get that $u(t)\in \LL^\infty(\Omega)$ for all $t>0$\,.

\noindent(iv) \textit{Nonnegative $H^*$ solutions are strong solutions. }In the case that the initial datum $u_0\in H^*$ is nonnegative, the $H^*$ solutions are bounded weak energy solutions, and by similar arguments to those given in Section 8.1 of \cite{DPQRV2} we can conclude that such solutions are strong in the sense of of Definition 4.1 of \cite{DPQRV2}, namely that $u_t\in \LL^\infty((\tau,\infty): \LL^1(\Omega))$ for all $\tau>0$\,.

\section{Asymptotic behaviour. Convergence to stationary state}\label{sec.asymp}

We begin the section by proving convergence in $\LL^\infty$ to the stationary profile and existence for the stationary equation, summarized in Theorems \ref{convLinf} and \ref{thStat}.
We begin with the proof of Theorem \ref{convLinf}. The proof is really similar to the one of Theorem 1.1 of \cite{V2} therefore we will just sketch it. We will split the proof in several steps.

\noindent$\bullet~$\textsc{Step 1. }\textit{A priori estimates. }We recall now the absolute bounds of \cite{BV2013}\,, valid for a (very weak) solution to the Dirichlet problem \eqref{FPME.Problem.PHI}, namely there exists a constant $K_1>0$ such that
\begin{equation}\label{Thm1.Step.1.1}
u(\tau,x)\le K_1\,\tau^{-\frac{1}{m-1}}\qquad\mbox{for a.e $(\tau,x)\in (0,\infty)\times \Omega$\,.}
\end{equation}
We also recall the monotonicity estimates
\begin{equation}\label{Thm1.Step.1.2}
u_\tau\ge -\frac{u}{(m-1)\tau}\qquad\mbox{for a.e. $(\tau,x)\in (0,\infty)\times \Omega$\,,}
\end{equation}
which have been first proved for the case $s=1$ in \cite{BCr-cont}\,, see also \cite{VazBook}\,. For the case $s\in (0,1]$\,, and more general linear operators and nonlinearities, see \cite{CP-JFA} and also \cite{DPQRV2, BV2013, BV2014}\,.  It is useful to pass to the rescaled solution
\[
v(t,x)=\tau^{1/(m-1)}u(\tau, x)\qquad\mbox{with}\qquad t=\log\tau\,.
\] In this way, the above estimates \eqref{Thm1.Step.1.1} and \eqref{Thm1.Step.1.2} transforms into
\begin{equation}\label{Thm1.Step.1.3}
0\le v(t,x)\le K_1\qquad\mbox{and}\qquad v_t\ge 0 \qquad\mbox{for a.e. $(t,x)\in (0,\infty)\times \Omega$\,.}
\end{equation}
Notice that $t\ge 0$ corresponds to $\tau\ge 1$\,.

\medskip

\noindent$\bullet~$\textsc{Step 2. }\textit{Convergence. }The estimates \eqref{Thm1.Step.1.3} valid for the solution $v$ to the time rescaled flow \ref{FPME.prob.log}\,, allow to conclude that for every $x\in \Omega$ there exists the limit
\[
\lim_{t\to\infty}v(t,x)=S(x)
\]
and this convergence is monotone nondecreasing, therefore the limit is nontrivial, $S\not\equiv 0$\,. Moreover, estimates \eqref{Thm1.Step.1.3} show that $S$ is bounded by $K_1$\,. Moreover, by Beppo-Levi Monotone Convergence Theorem we have also that
\[
v(t,\cdot)\to S \qquad\mbox{and}\qquad v^m(t,\cdot)\to S^m\qquad\mbox{in $\LL^1(\Omega)$}
\]
Therefore, since there is a uniform $\LL^\infty$ bound, the convergence takes place in any $\LL^p$\,, with $p<\infty$\,.

\noindent The $\LL^\infty$ convergence follows by $C^\alpha$ regularity, that in the case $s=1$ is well known, while for $0<s<1$ it has been proven in \cite{AC}, see also \cite{DPQRV2}. See also Subsection \ref{ssec.Calpha-Linfty}\,.

\noindent This concludes the proof of Theorem \ref{convLinf}. We now come to the proof of Theorem \ref{thStat}.

By our definition of $H^*$-solutions, we know that $v^m\in H^*(\Omega)$, so multiplying the equation \eqref{FPME.prob.log} by a test function $\psi\in H(\Omega)$, and integrating in $(t,t+T)\times\Omega$\,, with $T,t>0$\,, we obtain

\[\begin{split}
\int_\Omega v(t+T,x)\psi(x)\dx-\int_\Omega v(t,x)\psi(x)\dx
&=-\int_{t}^{t+T}\int_\Omega \A(v^m)\psi\dx\dt+\frac{1}{m-1}\int_{t}^{t+T}\int_\Omega v\psi\dx\dt\\
&=-\int_{t}^{t+T}\int_\Omega v^m \, \A\psi\dx\dt+\frac{1}{m-1}\int_{t}^{t+T}\int_\Omega v\psi\dx\dt
\end{split}
\]
where in the last step we have just used the integration by parts formula \eqref{Int.by.Parts}. Now we keep $T>0$ fixed and we let $t\to \infty$\,, so that the first member goes to zero\,, while the second member converges to
\[
-T\int_\Omega S^m \, \A\psi\dx\dt+\frac{T}{m-1}\int_\Omega S\psi\dx\dt
\]
Therefore we have that, in the limit $t\to \infty$
\[
\int_\Omega S^m \, \A\psi\dx\dt=\frac{T}{m-1}\int_\Omega S\psi\dx\dt
\]
which is the $H^*$ formulation of the solution to the elliptic problem \eqref{FPME.stat}\,.

Moreover, the stationary solution is unique: assume that we have two stationary solutions $S_1$ and $S_2$. Then we can construct two solutions of separation of variables for the non-rescaled flow as follows
\[
U_1(t,x)=\frac{S_1(x)}{t^{\frac{1}{m-1}}}\qquad\mbox{and}\qquad U_2(t,x)=\frac{S_2(x)}{(t+t_2)^{\frac{1}{m-1}}}
\]
for some $t_2>0$\,. Notice that $U_2$ has initial data $U_2(0,x)=S_2(x)\,t_2^{-\frac{1}{m-1}}\le U_1(0,x)=+\infty$\,. Therefore by comparison $U_2(t,x)\le U_1(t,x)$ for all $t>0$ and $x\in \Omega$\,, which can be expressed in the form
\[
S_2(x)\le \left(\frac{t+t_2}{t}\right)^{\frac{1}{m-1}}S_1(x)\qquad\mbox{for all $x\in \Omega$ and $t,t_2\ge 0$\,.}
\]
Finally, we let $t_2\to 0^+$ to get $S_2\le S_1$\,. The proof of the converse inequality $S_1\le S_2$ is obtained in a similar way.\qed

\section{The friendly giant and $\LL^\infty$-convergence with rates}\label{sec.fg}

In this section we prove the convergence with sharp rate in $\LL^\infty$, namely Theorem \ref{thm.PME.Rates.Giant}. The proof is based on the Global Harnack principle (GHP) proven in \cite{BV2013}, that we recall in the next subsection. The GHP holds for a suitable class of weak solutions, namely weak dual solutions, introduced in \cite{BV2013}, which a priori are different from the class of solutions we are dealing with. Indeed, we will show in the next subsection that weak dual solutions are indeed $H^*$-solutions.

 As a consequence of the GHP, we can construct two selfsimilar solutions of the form $U(t+\tau,x)=(t+\tau)^{-1/(m-1)}S(x)$\,, with $\tau\ge 0$. The upper bounds of the GHP, imply the existence of $U(t,x)=t^{-1/(m-1)}S(x)$ (case $\tau=0$)\,, which  corresponds to the initial datum $U(0,x)=+\infty$\,, and it has been called \textit{Friendly Giant}, cf. \cite{DahlKenig,V2, VazBook}, since it represents a upper barrier for any other solution. On the other hand, the lower bounds of the GHP, hold only after a waiting time $t_*=t_*(u_0)$\,, and imply that any solution stays above a solution of the form $U(t+c\,t_*,x)$ for some $c>1$\,. As a consequence, for large time, we show that any solution can be trapped between a Friendly Giant and lower selfsimilar solution obtained by time displacement. In this way  we can derive our asymptotic estimate, which turns out to be sharp, as one can observe by taking a time-displaced  solution of the form $U_1(t,x)=(t+t_1)^{-1/(m-1)}S(x)$ with $t_1>0$\,:
\[
\left|\frac{U_1(t,x)}{U(t,x)}-1\right|=1-\left(\frac{t}{t+t_1}\right)^{\frac{1}{m-1}}
=1- \left(1-\frac{t_1}{t+t_1}\right)^{\frac{1}{m-1}}
\le \frac{2\,t_1}{(m-1)(t_1+t)}
\]
$1-(1-y)^\beta\le 2\beta y$ valid for $y\in (0,1/2)$ and $\beta>0$ with the choice $y=t_1/(t+t_1)\le 1/2$, which corresponds to $t\ge t_1$\,.

\subsection{Review of quantitative estimates for the Dirichlet problem}\label{GHP.Estimates}

Let us recall the Global Harnack Principle, which encodes both absolute upper and lower bounds for the Dirichlet problem we are dealing with, and whose proof can be found in \cite{BV2013}. Before doing that, let us recall some definitions and setup some notations.

\noindent We denote the weighted $\LL^p$ space as $\LL^p_{\Phi_1}(\Omega)=\LL^p(\Omega\,,\, \Phi_1\dx)$, endowed with the norm
\[
\|f\|_{\LL^p_{\Phi_1}(\Omega)}=\left(\int_{\Omega} |f(x)|^p\Phi_1(x)\dx\right)^{\frac{1}{p}}\,.
\]
We will denote moreover by $\A$ the fractional operator which can be either the Spectral or the Restricted fractional Laplacian. We denote by $(\lambda_k, \Phi_k)$  its eigenvalues written in increasing order and repeated according to their multiplicity and its corresponding normalized eigenfunctions respectively. We also need a definition of suitable weak solutions introduced in \cite{BV2013}\,.

\begin{defn}[Weak dual solutions]
We say that $u \in C([0,T], L^1_{\Phi_1}(\Omega))$ is a weak dual solution if $\varphi(u) \in L^1([0,T], L^1_{\Phi_1}(\Omega))$ and
\begin{equation}
\int_0^T \int_\Omega \A^{-1}u \psi_t\dx\dt= \int_0^T \int_\Omega \varphi(u) \A \psi\dx\dt\,,\qquad\forall\,\psi/\Phi_1 \in C^1_c([0,T], L^\infty(\Omega))\,.
\end{equation}
\end{defn}
 \noindent This is a formulation of weak solutions for the potential equation $\partial_t \AI u+\varphi(u)=0$\,, which is a priori different from the $H^*$ solutions of Definition \ref{def.H*sols}. We will show in Lemma \ref{lem.5.3} that indeed $H^*$ solutions are weak dual solutions.

\begin{thm}[Global Harnack Principle for PME, \cite{BV2013}]\label{thm.GHP.PME}
Let let $m>1$ and let $u$ be a weak dual solution to the Dirichlet problem \eqref{FPME.Problem.PHI} corresponding to the initial datum $u_0\in \LL^1_{\Phi_1}(\Omega)$\,,  Then, there exists a constant $H_0\,,\,H_1>0$ and a time $t_*>0$ of the form
\begin{equation}\label{thm.GHP.PME.0}
t_*:=\tau_* \left[\frac{\int_{\Omega}\Phi_1\dx}{\int_{\Omega}u_0\Phi_1\dx }\right]^{m-1}\,,
\end{equation}
such that for all $t\ge t_*$ and $x_0\in \Omega$ the following inequality holds:
\begin{equation}\label{thm.GHP.PME.1}
H_0\,\frac{\Phi_1(x_0)^{\frac{1}{m}}}{t^{\frac{1}{m-1}}} \le \,u(t,x_0)\le H_1\, \frac{\Phi_1^{\frac{1}{m}}(x_0)}{t^{\frac{1}{m-1}}}
\end{equation}
The constants $\tau_*>0$ and $H_1,H_0>0$ depend only on $d, m, s$ and $\Omega$\,, but not on $u$\,. Recall that $\Phi_1$ is the first eigenfunction of $\A$ so that $k_{0,\Omega}\, \big(\dist(\cdot, \partial\Omega)^\gamma\wedge 1\big)\le \Phi_1 \le k_{1,\Omega}\, \big(\dist(\cdot, \partial\Omega)^\gamma\wedge 1\big)\,,$ with $\gamma=1$ for the SFL and $\gamma=s$ for the RFL\,.
\end{thm}

\noindent\textbf{Remark. }The previous estimates hold for weak dual solutions which a priori are a different class from the solutions at hand, namely $H^*$ solutions. We shall prove that $H^*$ solutions are indeed weak dual solutions.

\begin{lem}\label{lem.5.3}
Let $u$ be an $H^*$ solution. Then $u$ is a weak dual solution.
\end{lem}
\noindent {\sl Proof.~}Observe that the class of admissible test functions for $H^*$ solutions is bigger than the one of weak dual solutions. Indeed, if $\psi/\Phi_1 \in C^1_c((0,+\infty), L^\infty(\Omega))$ then $\psi \in C^1_c((0,+\infty), H^*(\Omega))$:
$$
\|\psi\|^2_{H^*}=\int_\Omega \psi \A^{-1} \psi\leq \|\psi/\Phi_1\|_\infty \int_\Omega \Phi_1 \A^{-1} \psi\leq\lambda_1^{-1} \|\psi/\Phi_1\|^2_\infty.
$$
Therefore, it is sufficient to show that
\begin{equation}\label{temp}
\int_0^T \int_\Omega u \partial_t (\A^{-1}  \psi)= \int_0^T \int_\Omega \A^{-1} u \partial_t \psi= \int_0^T \int_\Omega u \A^{-1}(\partial_t \psi)\,,\,\,\,\forall\, \psi/\Phi_1 \in C^1_c([0,T], L^\infty(\Omega))
\end{equation}
We have that
$$\|u\|_{L^1_{\Phi_1}} \leq \|u\|_{H^*}.$$
Indeed, we have
$$\|u\|^2_{L^1_{\Phi_1}}=\int_\Omega |u|\Phi_1=\lambda_1 \int_\Omega |u|(\A)^{-1}\Phi_1=\lambda_1\int_\Omega \A^{-1/2}u(\A)^{-1/2}\Phi_1. $$
By Cauchy-Schwarz inequality, this leads
$$\|u\|^2_{L^1_{\Phi_1}} \leq \lambda_1 \Big ( \int_\Omega \A^{-1/2}u\A^{-1/2}u\Big )^\frac12 \Big ( \int_\Omega \A^{-1/2}\Phi_1 \A^{-1/2}\Phi_1\Big )^\frac12=$$
$$\Big ( \int_\Omega u\A^{-1}u\Big )^\frac12\Big ( \lambda^2_1\int_\Omega \Phi_1 \A^{-1}\Phi_1\Big )^\frac12.$$
Hence the result. Furthermore, one has
$$\|\partial_t \psi\|^2_H=\int_\Omega \partial_t \psi \A(\partial_t \psi) \leq \|\psi/\Phi_1\|_\infty \int_\Omega \Phi_1 (\A(\partial_t \psi)\leq \lambda_1 \|\psi/\Phi_1\|^2_\infty .$$
Therefore, $\partial_t \psi \in H$ so \eqref{temp} holds. This proves the lemma. \qed

\subsection{The construction of the Friendly Giant. Proof of Theorem \ref{thm.PME.Rates.Giant}}

Now we are ready to prove Theorem \ref{thm.PME.Rates.Giant}. This will be a consequence of the following lemmata.

\begin{lem}[The Friendly Giant]
Under the assumptions of Theorem $\ref{thm.GHP.PME}$ we have that
\begin{equation}\label{Friendly.Giant.estimate}
u(t,x)\le U(t,x)=\frac{S(x)}{t^{\frac{1}{m-1}}}\,,\qquad\mbox{for all }(t,x)\in (0,\infty)\times\Omega\,.
\end{equation}
\end{lem}
\noindent {\sl Proof.~}Consider the solution by separation of variables:
\begin{equation}
U(\tau,x)=\frac{S(x)}{\tau^{\frac{1}{m-1}}}\,.
\end{equation}
This formally corresponds to the initial datum $u_0=+\infty$\,, therefore by comparison we know that it dominates all the other solutions, more precisely that inequality \eqref{Friendly.Giant.estimate} holds. This argument is formal, and we shall make it rigourous. Consider a sequence of initial datum $u_{0,n}=n$ and the corresponding weak dual solutions $u_n(t,x)$ which satisfy the absolute upper bounds \eqref{thm.GHP.PME.1} that are known to hold for all $t>0$, cf. \cite{BV2013}\,, namely
\[
u_n(t,x)\le \frac{K_1}{t^{\frac{1}{m-1}}}\,,\qquad\mbox{for all }(t,x)\in (0,\infty)\times\Omega\,.
\]
The family $u_n$ is monotone increasing by comparison, therefore the limit exists and it is finite thanks to the above absolute bounds, namely:
\[
U(t,x)=\lim_{n\to\infty}u_n(t,x)\le \frac{K_1}{t^{\frac{1}{m-1}}}\,.
\]
Moreover, if $u$ is a solution, then $\widetilde{u}(t,x)=k u(k^{m-1}t,x)$ is also a solution corresponding to $\widetilde{u}(0,x)=k u_0$\,. Hence the sequence $\widetilde{u_n}$ is monotone and converges to the same limit $U$, because $\widetilde{u_n}=u_{kn}$\,, and
\[
\lim_{n\to\infty}u_n(t,x)=\lim_{n\to\infty}\widetilde{u_n}(t,x)
\]
implies that $U$ is scaling invariant, i.e. $U(t,x)=k U(k^{m-1}t,x)$ for all $k>0$ and $(t,x)\in (0,\infty)\times\Omega\,.$ Finally we can put $U(t,x)=t^{-1/(m-1)}U(1,x)=t^{-1/(m-1)}S(x)$\,, where $S$ is the stationary solution.\qed

\begin{lem}
Under the assumptions of Theorem $\ref{thm.GHP.PME}$ we have that there exist
\begin{equation}\label{t0.t.star}
t_*:=\tau_* \left[\frac{\int_{\Omega}\Phi_1\dx}{\int_{\Omega}u_0\Phi_1\dx }\right]^{m-1}\qquad\mbox{and}\qquad
t_0=\k_0\, t_*\,,
\end{equation}
such that
\begin{equation}\label{Lower.Giant.estimate}
u(t,x)\ge U(t+t_0,x)=\frac{S(x)}{(t+t_0)^{\frac{1}{m-1}}}\,,\qquad\mbox{for all }(t,x)\in (t_*,\infty)\times\Omega\,.
\end{equation}
where $\k_0,\tau_*>0$ only depend on $d, s, m$\,.
\end{lem}
\noindent {\sl Proof.~}Recall that $S(x)$ is the solution to the stationary elliptic problem and satisfies the inequalities $h_0\Phi_1(x_0)^{\frac{1}{m}} \le S(x)\le h_1\Phi_1(x_0)^{\frac{1}{m}} $ for suitable constants $h_0, h_1>0$\,, cf. Theorem 9.2 of \cite{BV2013}\,, therefore the absolute lower bounds of Theorem \ref{thm.GHP.PME} at $t=t_*$ can be rewritten as follows:
\begin{equation}\label{lower.Giant}
u(t_*,x)\ge H_0\,\frac{\Phi_1(x_0)^{\frac{1}{m}}}{t_*^{\frac{1}{m-1}}}
\ge \frac{H_0}{h_1} \frac{S(x)}{t_*^{\frac{1}{m-1}}}
=\frac{S(x)}{(t_*+t_0)^{\frac{1}{m-1}}}=U(t_*+t_0,x)
\end{equation}
where in the last step we have chosen
\[
t_0=\left(  \frac{h_1^{m-1}}{H_0^{m-1}}  -1\right)t_*:= \k_0\, t_*\qquad\mbox{where}\qquad
t_*:=\tau_* \left[\frac{\int_{\Omega}\Phi_1\dx}{\int_{\Omega}u_0\Phi_1\dx }\right]^{m-1}\,.
\]
Therefore, by comparison, we have that $u(t,x)\ge U(t+t_0,x)$ for all $t\ge t_*$\,.\qed

\noindent {\bf Proof of Theorem \ref{thm.PME.Rates.Giant}.} The results of Lemma \ref{Friendly.Giant.estimate} and \ref{Lower.Giant.estimate} together, say that there exist a $t_0,t_*$ as in \eqref{t0.t.star} such that for all $t\ge t_*$  and $x\in \Omega$ we have
\[
\frac{S(x)}{(t+t_0)^{\frac{1}{m-1}}}\le u(t,x)\le \frac{S(x)}{t^{\frac{1}{m-1}}}\,,\qquad\mbox{or equivalently}\qquad
\left(\frac{t}{t+t_0}\right)^{\frac{1}{m-1}}\le \frac{u(t,x)}{U(t,x)}\le 1\,,
\]
which implies that
\[
\left|\frac{u(t,x)}{U(t,x)}-1 \right|\le 1- \left(\frac{t}{t+t_0}\right)^{\frac{1}{m-1}}
=1- \left(1-\frac{t_0}{t+t_0}\right)^{\frac{1}{m-1}}
\le \frac{2}{m-1}\frac{t_0}{t_0+t}
\]
where in the last step we have used the numerical inequality $1-(1-y)^\beta\le 2\beta y$ valid for $y\in (0,1/2)$ and $\beta>0$ with the choice $y=t_0/(t-t_0)\le 1/2$, which corresponds to $t\ge t_0$\,. Therefore we have proved the desired bound for all $t\ge t_0$\,, where $t_0$ has the  expression
\[
t_0=\left(  \frac{h_1^{m-1}}{H_0^{m-1}}  -1\right)t_*
=\left(  \frac{h_1^{m-1}}{H_0^{m-1}}  -1\right)\tau_* \left[\frac{\int_{\Omega}\Phi_1\dx}{\int_{\Omega}u_0\Phi_1\dx }\right]^{m-1}
=\overline{k}\left[\frac{\int_{\Omega}\Phi_1\dx}{\int_{\Omega}u_0\Phi_1\dx }\right]^{m-1}
\]
where the constants $H_0, h_1, \tau_*>0$ have an explicit expression given in Theorems 8.1 and 9.2 of \cite{BV2013}\,.\qed

\section{Entropy dissipation and relative error: proof of Theorem \ref{PME.Rates}}\label{sec.entr}

We now would like to pass to the relative error $w(t,x)=\frac{v(t,x)}{S(x)}-1$, but it turns out to be impossible to write a clean equation due to the nonlocality of $\A$, without passing through the extension technique by Caffarelli-Silvestre \cite{CaSi} in the elliptic case, see also \cite{CabTan, c-d-d-s,RosSer}\,; as for the extension technique in the parabolic case, see \cite{AC} and \cite{DPQRV1,DPQRV2}.

\subsection{The relative error and its extension}

One has to consider two different extensions: one for the parabolic problem and one for the stationary problem. Furthermore, the setting is slightly different if one considers the RFL or the SFL.

\subsubsection*{The extension property for the SFL}

It has been proved in \cite{c-d-d-s} that the following holds.

\begin{lem}
Consider a weak solution of
\begin{equation}
\left\{
\begin{array}{lll}
\mbox{\rm div}(y^{1-2s} \nabla w)=0 &\mbox{in }\,\,\mathcal C= \Omega \times (0,+\infty),\\
w=0\,,\; &\mbox{on }\,\,\partial \Omega \times (0,+\infty) \\
\end{array}
\right.
\end{equation}
Then $-\lim_{y \to 0} y^{1-2s}\partial_y w=\mathcal L w(\cdot,0).$
where $\mathcal L$ is the SFL.
\end{lem}

\subsubsection*{The extension property for the RFL}

It has been proved in \cite{CaSi} that the following holds.

\begin{lem}
Consider a weak solution of
\begin{equation}
\left\{
\begin{array}{lll}
\mbox{\rm div}(y^{1-2s} \nabla w)=0 &\mbox{in }\,\, \mathbb R^n \times (0,+\infty),\\
w=0\,,\; &\mbox{on }\,\,(\mathbb R^n\setminus\Omega)\times\{0\} \\
\end{array}
\right.
\end{equation}
Then $-\lim_{y \to 0} y^{1-2s}\partial_y w=\mathcal L w(\cdot,0).$
where $\mathcal L$ is the RFL.
\end{lem}

While the elliptic extension above has been widely used, for the parabolic we have to be a bit more careful. Let us notice first that due to the bounds in Theorem \ref{thm.GHP.PME}, any $H^*$ solution is uniformly bounded for all $t>0$. Hence, in particular, it is a weak solution in the sense of \cite{DPQRV2} in the case of the SFL. Indeed the authors of \cite{DPQRV2} have proved that any weak solution of the parabolic problem is equivalent to a weak solution of the Caffarelli-Silvestre associated extension. Similar arguments hold also in the case of the RFL. This fully justifies the computations below. We will treat both operators at the same time.

Denote by $\tilde \Omega$ either the cylinder $\mathcal C$ for the SFL or the whole half-space $\mathbb R^{n+1}_+$ for the RFL. We first extend the solution $v^m$ to a function $\hat{V}(t,x,y)$ defined for $(x,y)$ in the upper cylinder $Q= \tilde \Omega\times \RR_+$ and $t>0$: ($\hat V=\mbox{Ext} (v^m) )$
\begin{equation}\label{FPME.prob.log.ext}
\begin{split}
\left\{\begin{array}{lll}
\text{div}(y^{1-2s} \nabla \hat{V})=0 & ~ {\rm in}~\RR_+\times Q,\\[3mm]
\partial_tv-y^{1-2s}\partial_y \hat{V} =c\hat V^{1/m} & ~ {\rm in}~\RR_+\times \Omega,\\[3mm]
\end{array}\right.
\end{split}
\end{equation}
Analogously we extend the solution to the stationary problem $S^{m}$ to a function $\hat S$ on the upper cylinder $Q$:( $\hat S=Ext(S^m)$)
\begin{equation}\label{FPME.stat.ext}
\begin{split}
\left\{\begin{array}{lll}
\text{div}(y^{1-2s} \nabla \hat{S})=0 & ~ {\rm in} \ Q,\\[3mm]
-y^{1-2s}\partial_y \hat S =\dfrac{\hat S^{1/m}}{m-1} & ~ {\rm in}~\Omega,\\[3mm]
\end{array}\right.
\end{split}
\end{equation}

In the previous problems, depending which operator $\A$ we are considering, the boundary conditions are the following:
\begin{itemize}
\item $\hat V = \hat S =0 $ on $\mathbb R^+  \times \partial \Omega$ in the case of  SFL
\item $\hat V = \hat S =0 $ on $( \partial \mathbb R^{N+1}_+ \backslash (\mathbb R^+ \times \Omega))$ in the case of the RFL
\end{itemize}
Furthermore, the following generalized integration by parts formula holds
\begin{equation}\label{IPP}
\int_{\tilde \Omega} \psi_1 y^{1-2s} \partial_y \psi_2= -\int_{\tilde \Omega \times \mathbb R^+} \psi_1\nabla \cdot (y^{1-2s} \nabla \psi_2)-\int_{\tilde \Omega \times \mathbb R^+} y^{1-2s} \nabla \psi_1 \cdot \nabla \psi_2.
\end{equation}
Now we are ready to define the \textit{extended relative error} on the upper cylinder $Q$ as
\begin{equation}\label{ext.rel.err}
W(t,x,y)=\left(\frac{\hat{V}(t,x,y)}{\hat S(x,y)}\right)^{\frac{1}{m}}-1\qquad\mbox{i.\,e.}\qquad \hat{V}=\hat S(W+1)^m\,.
\end{equation}
We now compute an equation for the extended relative error $W$. We have
$$\text{div}(y^{1-2s}\nabla \hat V)=0= m \, \text{div}(y^{1-2s}\hat S (1+W)^{m-1}\nabla W)+\text{div}(y^{1-2s}(1+W)^m \nabla \hat S). $$
We now use the general formula
$$\text{div}(\eta\Upsilon)=\nabla \eta\cdot \Upsilon+( \text{div} \Upsilon )\,\eta$$
where $\eta$ is a real-valued function and $\Upsilon$, a vector-valued one. We obtain
$$0=m\, \text{div}(y^{1-2s}\nabla W) (1+W)^{m-1} \hat S+ m\, y^{1-2s}(1+W)^{m-1}\nabla W\cdot \nabla \hat S+$$
$$m(m-1)\, y^{1-2s}(1+W)^{m-2}|\nabla W|^2 \hat S + m\, y^{1-2s}(1+W)^{m-1}\nabla W\cdot \nabla \hat S . $$
This  shows that the function $W$ satisfies the equation
\begin{equation}
  \text{div}(y^{1-2s} \nabla W)=-y^{1-2s}(m-1)\dfrac{\left|\nabla W\right|^2}{W+1}-2y^{1-2s} \frac{\nabla \hat S\cdot\nabla W}{\hat S}\,.
\end{equation}
This is supplemented by a dynamical boundary condition in terms of the trace of the extended relative error, namely on the function $w=\tr(W)$:
\begin{equation}\label{wt}
\partial_t w= \frac{F(w)}{m-1}+mS^{m-1}(1+w)^{m-1}y^{1-2s}\partial_yW|_{y=0}
\end{equation}
where we have introduced the function
\begin{equation}
F(w):=\left[(w+1)-(w+1)^m\right]\,.
\end{equation}

\subsection{The entropy dissipation}
Define the following entropy functional on functions $w\in\LL^2(\Omega)$
\begin{equation}\label{Entropy.L2}
\mathcal{E}[w]:=\frac{1}{2}\int_\Omega\left|w(x)-\w\right|^2\,S^{1+m}\dx\qquad\mbox{with}\qquad \w:=\frac{\int_\Omega w\,S^{1+m}\dx}{\int_\Omega S^{1+m}\dx}\,.
\end{equation}
Now we want to apply this functional to $w(t,x)=\tr(W(t,x,y))$, where $W$ is the extended relative error defined in \eqref{ext.rel.err}, and calculate the time derivative along the flow, which is the entropy dissipation, and it will be related to the extended relative error $W$. Recall that $\tilde \Omega$ is either the cylinder $\mathcal C$ for the SFL or the whole half-space $\mathbb R^{n+1}_+$ for the RFL, and $Q= \tilde \Omega\times \RR_+$.
\begin{prop}[Entropy production]\label{entropy.1} Under the running assumption, we have that for all $t>0$
\begin{equation}\label{entr.prod}
\frac{\rd}{\dt}\mathcal{E}[w](t)
    =-m\int_{Q} y^{1-2s} (W+1)^{m-1}\left|\nabla W\right|^2 \, \hat S^2\dx\dy\,
    +\frac{1}{m-1}\int_\Omega\big(w(x)-\w\big)F(w)S^{m+1}\dx\,.
\end{equation}
\end{prop}

\noindent {\sl Proof.~} All the following computations are justified in view of the fact that solutions $u$ at our disposal are weak energy solutions and they are also strong (since nonnegative), as explained in Remarks (iii) and (iv) of Section \ref{sec-Brezis}. Therefore also the rescaled  solutions $v$ enjoy similar regularity properties and also the relative error $w=(v/S)^m-1$\,.   We have
$$
\frac{\rd}{\dt}\mathcal{E}[w](t)
    =\frac{1}{2}\frac{\rd}{\dt}\int_\Omega\left|w(x)-\w\right|^2\,S^{1+m}\dx=
    \int_\Omega\big(w(x)-\w\big)\,S^{1+m}\, ( \partial_t w(x)-\partial_t \w(x))\dx
$$
Let us put $\rd\mu=S^{1+m}\dx$. We remark that
\begin{equation}
\begin{split}
\int_\Omega\big[w(t)-\w(t)\big]\,\big[\partial_t\w(t)\big]\,\rd\mu
& =\big[\partial_t\w(t)\big]\int_\Omega\big[w(t)-\w(t)\big]\,\rd\mu =0\,,
\end{split}
\end{equation}
and also that
\begin{equation}\label{var}
0\le \int_\Omega \big|w-\w(t)\big|^2 \rd\mu
    =\int_\Omega \big[w^2-\w(t)^2\big]\rd\mu\,.
\end{equation}
Using also Equation \eqref{wt} we get
\[
\frac{\rd}{\dt}\mathcal{E}[w](t)=
    m\int_\Omega\big(w(x)-\w\big)\,S^{2m}\,(w+1)^{m-1}y^{1-2s}\partial_y W \dx
        + \frac{1}{m-1}\int_\Omega\big(w(x)-\w\big)F(w)S^{m+1}\dx\,.
\]
Next we calculate the first of the integrals in the RHS:
\[
\begin{split}
&\int_\Omega\big(w(x)-\w\big)\,S^{2m}\,(w+1)^{m-1} y^{1-2s}\partial_y W|_{y=0} \dx\\
 (a)&=-\int_{Q}y^{1-2s}\nabla \left[\big(W(x,y)-\w\big)\,\hat S^2(x,y)\,\big(W(x,y)+1\big)^{m-1}\right]\cdot\nabla W \dx\dy\,\\
    & -\int_{Q}\big(W(x,y)-\w\big)\,\hat S^2(x,y)\,\big(W(x,y)+1\big)^{m-1} \text{div}(y^{1-2s}  \nabla W(x,y)) \dx\dy\\
    &=-\int_{Q}y^{1-2s}\left|\nabla W\right|^2 \,(W+1)^{m-1} \hat S^2\dx\dy\,\\
    &-2\int_{Q}y^{1-2s}(W-\w\big)\,(W+1)^{m-1}\hat S\nabla \hat S \cdot\nabla W \dx\dy\, \\
    &-(m-1)\int_{Q}y^{1-2s} (W-\w\big)\,(W+1)^{m-2}\left|\nabla W\right|^2 \, \hat S^2\,   \dx\dy\, \\
    \end{split}
\]
\[
\begin{split}
 (b)& +(m-1)\int_{Q} y^{1-2s}(W-\w\big)\,\hat S^2\,(W+1)^{m-1}\dfrac{\left|\nabla (W+1)\right|^2}{W+1}\dx\dy \\
    & +2\int_{Q}y^{1-2s} (W-\w\big)\,\hat S\,(W+1)^{m-1}\nabla\hat S\cdot\nabla W\dx\dy \\
    &=-\int_{Q}y^{1-2s} (W+1)^{m-1}\left|\nabla W\right|^2 \, \hat S^2\dx\dy\,\\
\end{split}
\]
where in (a) we have integrated by parts \eqref{IPP} and in $(b)$ we have used the equation satisfied by $W$ on $Q$. Summing up we have obtained formula \eqref{entr.prod}\,.\qed

First we estimate the second integral of \eqref{entr.prod}: we just need the trace $w=\tr(W)$ of the extended relative error $W$ to be small, and this is true for large times, as a consequence of Theorem \ref{convLinf}.

Let us fix $\varepsilon>0$ and define $t_0$ as the time for which .

\begin{lem}\label{entropy.2}
Let $m>0$. Under the running assumption we have that
\[
\|w(t)\|_{\LL^\infty(\Omega)}\le \varepsilon\qquad\mbox{for all $t\ge t_\varepsilon$}\,,
\]
and
\[
\begin{split}
\int_\Omega\big(w(t,x)-\w\big)F(w)S^{m+1}\dx
&\le -2(m-1)\big[1-\varepsilon(m+1)\big]\mathcal{E}[w](t)
\end{split}
\]
where $\w=\int_\Omega w\rd\mu$ and $F(w)=(w+1)-(w+1)^m$\,.
\end{lem}
\noindent {\sl Proof.~}Let us put $\rd\mu=S^{m+1}\dx$. We first notice that since $|w|\le\varepsilon$, then $|\w|\le \varepsilon$, and $|w-\w|\le 2\varepsilon$, so that Taylor expansion for $F$ at the point $\w$ gives
\[
F(w)=F(\w)+F'(\w)(w-\w)+\frac{1}{2}F''(\tilde{w})(w-\w)^2
\]
where $\tilde{w}$ lies between $w$ and $\w$. Hence we get
\[\begin{split}
\int_\Omega\big(w(t,x)-\w\big)F(w)\rd\mu
&= F(\w)\int_\Omega\big(w-\w\big)\rd\mu +F'(\w)\int_\Omega\big(w-\w\big)^2\rd\mu
+\frac{1}{2}\int_\Omega F''(\tilde{w})\big(w-\w\big)^3\rd\mu\\
&=F'(\w)\int_\Omega\big(w-\w\big)^2\rd\mu
+\frac{1}{2}\int_\Omega F''(\tilde{w})\big(w-\w\big)^3\rd\mu
\end{split}\]
Now we use the fact that
\[
F'(\w)=1-m(1+\w)^{m-1}=1-m + m(1-m)\w+ O(\w^2)\le 1-m + \varepsilon \big[m(1-m)+O(\varepsilon) \big]
\]
which tends to $(1-m)$ as $t\to \infty$, since we already know that $w(t)\to 0$ uniformly in space as $t\to \infty$, by Theorem \ref{convLinf}, hence $\w(t)\to 0$ as $t\to \infty$ as well. Also
\[
F''(\tilde{w}) = m(1-m)(1+\tilde{w})^{m-2}\to m(1-m)\qquad\mbox{and}\qquad F''(\tilde{w})=m(1-m)+O(\varepsilon)
\]
uniformly in space as $t\to\infty$ for the same reason. Putting these things together, we
have obtained that for any $m>0$
\[\begin{split}
\int_\Omega\big(w(t,x)-\w\big)F(w)\rd\mu
    &\le\Big[1-m+\varepsilon\left[m(1-m)+O\big((m-1)\varepsilon\big)\right]\Big]\, \int_\Omega\big(w(t,x)-\w\big)^2\rd\mu\\
    &\le-(m-1)\big[1-\varepsilon(m+1)\big]\int_\Omega\big(w(t,x)-\w\big)^2\rd\mu
\end{split}
\]
where in the last step we have chosen $\varepsilon$ so small that
\[
\Big[1-m+\varepsilon\left[m(1-m)+O\big((m-1)\varepsilon\big)\right]\Big]\le -(m-1)\big[1-\varepsilon(m+1)\big]
\]
which amounts to require that $t_0\gg 1$. Notice that in the limit $m\to 1$, the last term disappears, and this is exactly what happens in the linear case.\qed

\begin{prop}[Entropy decay rates]\label{entropy.3}

Let $m>1$. Under the running assumption we have that there exists a constant $K>0$ such that for all $t\ge t_0$
\[
\mathcal E[w](t) \leq K\, e^{-2(m-1)t}.
\]
where the constant $K$ depends on $m, \mathcal{E}[w](t_0)\,,$ but not on $w(t)$\,.
\end{prop}

\noindent {\sl Proof.~}From proposition \ref{entropy.1} and Lemma \ref{entropy.2} we get

\begin{equation}\label{entr.prod.2}
\frac{\rd}{\dt}\mathcal{E}[w](t)\le \frac{1}{m-1}\int_\Omega\big(w(x)-\w\big)F(w)S^{m+1}\dx
\le-2(m-1)\big[1-\varepsilon(m+1)\big]\mathcal{E}[w]
\end{equation}
then we integrate the above equation to get
\begin{equation}\label{entr.prod.3}
\mathcal{E}[w](t)\le \ee^{-2(m-1)\big[1-\varepsilon(m+1)\big](t-t_0)}\mathcal{E}[w(t_0)]
\end{equation}
Plugging \eqref{entr.prod.3} into \eqref{entr.prod.2} gives
\[
\frac{d\mathcal E [w](t)}{dt} \leq -2(m-1) \mathcal E[w](t) +c_1 \,\ee^{-c_2 (t-t_0)}
\]
where $c_1=2\varepsilon(m-1)(m+1)\mathcal{E}[w(t_0)]>0$ and $c_2=-2(m-1)\big[1-\varepsilon(m+1)\big]>0$.
The result comes from an integration.\qed

\noindent\textbf{Remark. }The entropy decay rate does not give directly information about the decay of $\LL^p$ norms.
We first need to check the decay of the weighted $\LL^1$-norm.

\begin{lem}[Weighted $\LL^1$-norm decay]\label{mean}
Let $m>1$. Under the running assumption we have that for all $t\ge t_0$
$$\overline w(t) \leq \ee^{-(m-1)t}\overline w_0. $$
\end{lem}
\noindent {\sl Proof.~}
We have
$$\partial_t \overline w(t)=\frac{1}{\int_\Omega S^{1+m}} \int_\Omega \partial_t w S^{1+m}\,dx. $$

\noindent Using Equation \eqref{wt}, we have

$$\partial_t \overline w(t)=\frac{1}{\int_\Omega S^{1+m}} \int_\Omega F(w) S^{1+m}\,dx +\frac{m}{\int_\Omega S^{1+m}}\int_{\Omega} (1+w)^{m-1} S^{2m} y^{1-2s} \partial_y W|_{y=0}\,dx .   $$

\noindent We now evaluate the integral $\int_{\Omega} (1+w)^{m-1} S^{2m} y^{1-2s} \partial_y W|_{y=0}\,dx . $ Formula \eqref{IPP}, gives
$$\int_{\Omega} (1+w)^{m-1} S^{2m} y^{1-2s} \partial_y W|_{y=0}\,dx=\int_{Q} (1+W)^{m-1}\hat S^2 \nabla \cdot (y^{1-2s} \nabla W)+\int_{Q} y^{1-2s} \nabla W\cdot \nabla (\hat S^2 (1+W)^{m-1}). $$
Using Equation \eqref{FPME.prob.log.ext}, we have
$$\int_{\Omega} (1+w)^{m-1} S^{2m} y^{1-2s} \partial_y W|_{y=0}\,dx=0.$$
So that
\[
\partial_t \overline w(t)=\frac{1}{\int_\Omega S^{1+m}} \int_\Omega F(w) S^{1+m}\,dx
\le -\frac{m-1}{\int_\Omega S^{1+m}}\int_\Omega F(w) S^{1+m}\,dx =-(m-1)\overline{w}\,,
\]
where we have used the concavity of $F$, namely $F(w)\le -(m-1)w$, as in \cite{BGV-Domains}. The result follows by integration.\qed

\begin{prop}[Weighted $\LL^2$-norm decay]\label{exp.decay.norm.PME}
Let $m>1$. Under the running assumption we have that there exists a constant $K_1>0$ such that for all $t\ge t_0$
\begin{equation}\label{rate.entropy.PME}
\int_\Omega\big|w(t)\big|^2 S^{1+m}\dx
\le K_1\ee^{-2(m-1)t}
\end{equation}
where the constant $K_1$ depends on $m, \varepsilon, \mathcal{E}[w(t_0)], \overline{w}_0 $ but not on $w(t)$\,.
\end{prop}
\noindent {\sl Proof.~}The result of Proposition \ref{entropy.3} can be rewritten as:
\begin{equation*}\begin{split}
\frac{1}{2}\int_\Omega\left|w(x)\right|^2\,S^{1+m}\dx&\le K\, \ee^{-2(m-1)t}+\frac{1}{2}\w^2\int_\Omega S^{1+m}\dx\\
&\le K\, \ee^{-2(m-1)t}+\frac{1}{2}\left(\int_\Omega S^{1+m}\dx\right)\ee^{-2(m-1)t}\overline w_0
\le K_1\ee^{-2(m-1)t}
\end{split}
\end{equation*}
where $K>0$ is as in Proposition \ref{entropy.3}\,, since we recall that
\[
\mathcal{E}[w]:=\frac{1}{2}\int_\Omega\left|w(x)-\w\right|^2\,S^{1+m}\dx
=\frac{1}{2}\int_\Omega\left|w(x)\right|^2\,S^{1+m}\dx - \frac{1}{2}\w^2\int_\Omega S^{1+m}\dx\,.
\]
and that $\overline w(t) \leq \ee^{-(m-1)t}\overline w_0, $ by Lemma \ref{mean}\,.\qed

\noindent\textbf{Remark. }The results of Lemma \ref{entropy.2} and Proposition \ref{entropy.3} have the following consequences.

\begin{cor}[Rates of decay for the difference]\label{L1.norm.PME}
Let $m>1$. Under the running assumption we have that there exists a constants $K_1,K_2>0$ such that for all $t\ge t_0$:
\begin{equation}\label{L1-weight}
\int_\Omega|v(t,\cdot)-S|\,S^{m} \dx \le \ee^{-(m-1)t}\int_\Omega|v_0-S|\,S^{m} \dx \,,
\end{equation}
\begin{equation}\label{L2-weight}
\int_\Omega|v(t,\cdot)-S|^2\,S^{m-1} \dx \le K_1\, \ee^{-2(m-1)t}\,,
\end{equation}
and
\begin{equation}\label{L1-norm}
\int_\Omega|v(t,\cdot)-S|\dx \le K_2\, \ee^{-2(m-1)t}\,.
\end{equation}
where the constant $K_1,K_2$ depends on $m, \gamma, u_0, v(t_0) $ but not on $v(t)$\,.
\end{cor}
\noindent {\sl Proof.~}Let recall that the results of Lemma \ref{entropy.2} and Proposition \ref{entropy.3} can be rewritten in the form \eqref{L1-weight} and \eqref{L2-weight} respectively. It only remains to prove \eqref{L2-weight}. Recall that by H\"older inequality, we have $\|f\|_{\LL^1(\Omega)}^2\le \|g\|_{\LL^1(\Omega)}\,\|f^2/g\|_{\LL^1(\Omega)}$\,, so that
\[\begin{split}
\left(\int_\Omega|v(t,\cdot)-S|\,\dx\right)^2\le \left(\int_\Omega|v(t,\cdot)-S|^2\,S^{m-1} \dx\right)\left(\int_\Omega S^{-(m-1)} \dx\right)
\le  K_2\, \ee^{-2(m-1)t}
\end{split}
\]
where we have used \eqref{L2-weight} and
\[
\int_\Omega S^{-(m-1)} \dx\le \frac{1}{h_0}\int_\Omega \Phi_1^{-\frac{m-1}{m}} \dx \le
\frac{k_{1,\Omega}}{h_0}\int_\Omega \big(\dist(x, \partial\Omega)^{-\gamma\frac{m-1}{m}} \dx := \frac{K_2}{K_1}\,,
\]
recall that $S(x)$ is the solution to the stationary elliptic problem and satisfies the inequalities $h_0\Phi_1(x_0)^{\frac{1}{m}} \le S(x)\le h_1\Phi_1(x_0)^{\frac{1}{m}} $ for suitable constants $h_0, h_1>0$\,, cf. Theorem 9.2 of \cite{BV2013}\,, therefore since $\Phi_1(x)\asymp\big(\dist(x, \partial\Omega)^\gamma\wedge 1\big)$\,, with $\gamma\in (0,1)$\,.\qed

\medskip

\subsection{$C^\alpha$ continuity, interpolation and rates}\label{ssec.Calpha-Linfty}

We first recall a useful interpolation Lemma due to Gagliardo \cite{Ga}, cf. also  Nirenberg, \cite[p. 126]{MR0109940}.  We leave the complete version for the reader's convenience.

\begin{lem}\label{GN} Let $\lambda$, $\mu$ and $\nu$ be such that
$-\infty<\lambda \le \mu
\le \nu <\infty$. Then there exists a positive constant $\mathcal
C_{\lambda,\mu,\nu}$ independent of $f$ such that
\begin{equation}\label{eq:interpolation}
\|f\|_{1/\mu}^{\nu-\lambda} \le \mathcal
C_{\lambda,\mu,\nu}\|f\|_{1/\lambda}^{\nu-\mu} \;
\|f\|_{1/\nu}^{\mu-\lambda}\,, \qquad\forall\;f\in C_c^{\infty}(\RR^d)\;,
\end{equation}
where $\|\cdot\|_{1/\sigma}$ stands for the following
quantities:\begin{itemize}
\item[(i)] If $\sigma>0$, then
$\|f\|_{1/\sigma}=\left(\int_{\RR^d}|f|^{1/\sigma}\dx\right)^\sigma$.

\item[(ii)] If $\sigma<0$, let $k$ be the integer part of
$(-\sigma d)$ and $\theta=|\sigma|d-k$ be the fractional
(positive) part of $\sigma$. Using the standard multi-index
notation, where $|\eta|=\eta_1+\ldots+\eta_d$ is the length
of the multi-index
$\eta=(\eta_1,\ldots\eta_d)\in\mathbb{Z}^d$, we define
\begin{equation*}\label{def.C^k}
\|f\|_{1/\sigma}=\left\{\begin{array}{lll}
\displaystyle\max_{|\eta|=k}\; \big|\partial^\eta
f\big|_\alpha=\displaystyle\max_{|\eta|=k}\;
\sup_{x,y\in\RR^d}\;\dfrac{\big|\partial^\eta f(x)-\partial^\eta
f(y)\big|}{|x-y|^\theta}&
\mbox{if~}\alpha>0\;,\\[5mm]
\displaystyle\max_{|\eta|=k}\;\displaystyle\sup_{z\in\RR^d}\big|\partial^\eta
f(z)\big|:=\|f\|_{C^k(\RR^d)}& \mbox{if~}\alpha=0\;.
\end{array}\right.
\end{equation*}
As a special case, we observe that $\sigma=-\alpha/d$ we get $\|f\|_{-d/\alpha}=\|f\|_{C^{\alpha}(\RR^d)}$ for any $\alpha\in (0,1]$.

\item[(iii)] If $\sigma=0$, by convention, we let
$\|f\|_{1/0}=\sup_{z\in\RR^d}|f(z)|=\|f\|_{C^0(\RR^d)}=\|f\|_{\infty}$.
\end{itemize}
\end{lem}

Next we show how the above result can be combined with $C^\alpha$ regularity in order to obtain the same rate of decay for all $\LL^p$-norms:

\begin{cor}
Let $v$ be the solution of $\ref{FPME.prob.log}$ and $S$ be the stationary solution. Assume that for some $t\ge \underline{t}>0$ the following estimates hold,
\begin{equation}\label{Calpha-Linfty}
\|v(t)-S\|_{C^\alpha(\Omega)}\le C\|v(t)-S\|_{\LL^\infty(\Omega)}\,.
\end{equation}
 for some constant $C>0$ which do not depend on $v(t)$.
Then there exists a constant $\overline{K}>0$ that does not depend on $u$ nor $S$ such that
\begin{equation}\label{cor.interp.Calpha}
\|u(t)-S\|_{\LL^\infty(\Omega)}\le \overline{K} \|u(t)-S\|_{\LL^1(\Omega)}\qquad\mbox{for all $t>t_*$.}
\end{equation}
\end{cor}
\noindent {\sl Proof.~}We use the interpolation inequality \eqref{eq:interpolation} with the choices: $\mu=0$\,, $\lambda=-d/\alpha$ and $\nu=1$\,, where $\alpha\in (0,1)$ is the regularity exponent of \eqref{Calpha-Linfty}, so that:
\begin{equation}\label{interp.Calpha}
\|u(t)-S\|_{\LL^\infty(\Omega)}^{\nu-\lambda}
\le \mathcal C_{\lambda,\mu,\nu}\|u(t)-S\|_{C^\alpha(\Omega)}^{\nu-\mu}\|u(t)-S\|_{\LL^1(\Omega)}^{\mu-\lambda}
\le \mathcal C_{\lambda,\mu,\nu} C_\alpha \|u(t)-S\|_{\LL^\infty(\Omega)}^{\nu-\mu}\|u(t)-S\|_{\LL^1(\Omega)}^{\mu-\lambda}
\end{equation}
where in the last step we have used the regularity estimate $\|u(t)-S\|_{C^\alpha(\Omega)}\le C_\alpha \|u(t)-S\|_{\LL^\infty(\Omega)}$\,, valid for $t\ge t_*$\,. The above inequality implies inequality \eqref{cor.interp.Calpha} with $C=\mathcal C_{\lambda,\mu,\nu} C_\alpha$.\qed

\noindent \textbf{Remarks. }\\(i) The above corollary allows to extend the sharp rates for the $\LL^1$-norm \eqref{L1-norm} to the $\LL^\infty$ norm.

\noindent  (ii) The $C^\alpha$ regularity expressed in the form of inequality \eqref{Calpha-Linfty}\,, is known to hold when $s=1$\,, see e.g. \cite{BGV-Domains}, while when $s\in (0,1)$ it may be derived as a consequence of the results of \cite{AC}\,.

\noindent (iii) If solutions are just $C^\alpha$, i.e. $\|u(t)-S\|_{\LL^\infty(\Omega)}\le C$\,, then inequality \eqref{interp.Calpha} implies that we can have a rate also for the $\LL^\infty$-norm, but we loose optimality of the exponents:
\[
\|u(t)-S\|_{\LL^\infty(\Omega)}\le \overline{K}\|u(t)-S\|_{\LL^1(\Omega)}^{\frac{d}{d+\alpha}}\,.
\]

%%%%%%%%%%%%%%%%%%%%%%%%%%%%%%%%%%%%%%%%%%%%%%%%%%%%%%%%%%%%%%%%%%%%%%%%%%%%
\section{Appendix 1. Reminder about fractional Sobolev spaces}

Throughout this section, we always consider a bounded Lipschitz domain $\Omega\subset \RR^d$, so that it enjoys the extension property. We will always consider $s\in(0,1)$\,, unless explicitly stated. When we write $ f\asymp g$ we mean that there exist two constants $0<c_0\le c_1$ such that $c_0 f\le g\le c_1 g$\,.

\subsection{The spaces $W^{s,2}(\RR^d)$ and $H^s(\RR^d)$. }
We first recall the definitions and basic properties of the fractional Sobolev spaces on the whole space: indeed in literature one can find several definitions which are well known to be all equivalent\,. The first definition is given in terms of the so-called Gagliardo seminorm:
\begin{equation}\label{Seminorm.Gagliardo}
[f]_{W^{s,p}(\RR^d)}:=\left(\int_{\RR^d}\int_{\RR^d} \frac{|f(x)-f(y)|^p}{|x-y|^{d+ps}}\dx\dy\right)^\frac{1}{p}\,.
\end{equation}
Then we define the space
\begin{equation}\label{Wps.Rd}
W^{s,p}(\RR^d):=\left\{f\in\LL^p(\RR^d)\,\big|\,[f]_{W^{s,p}(\Omega)}<+\infty \right\}
\end{equation}
which is a Banach space (Hilbert space for $p=2$) with the norm
\begin{equation}\label{Wps.Norm}
\|f\|_{W^{s,p}(\RR^d)}=  [f]_{W^{s,p}(\RR^d)}+\|f\|_{\LL^p(\RR^d)}
\end{equation}
Another possible definition, when $p=2$ is via the Fourier transform
\begin{equation}\label{Hs.Rd}
H^s(\RR^d):=\left\{f\in\LL^2(\RR^d)\,\big|\, \left(1+|\xi|^2\right)^{\frac{s}{2}}\hat{f}\in \LL^2(\RR^d)\right\}
\end{equation}
which is a Hilbert space, with the norm
\begin{equation}\label{Hs.Norm.Rd}
\|f\|_{H^s(\RR^d)}^2= \left\|\left(1+|\xi|^{2s}\right)^{\frac{1}{2}}\hat{f}\right\|_{\LL^2(\RR^d)}^2
\asymp \left\|\left(1+|\xi|^2\right)^{\frac{s}{2}}\hat{f}\right\|_{\LL^2(\RR^d)}^2
\end{equation}
See the book \cite{LM} pg. 35, where one can find also the proof that $H^s$ is an interpolation space between $H^1$ and $\LL^2$\,, namely $H^s(\RR^d)=\big[H^1(\RR^d)\,,\,\LL^2(\RR^d)\big]_{1-s}$ (Thm 7.1 of \cite{LM})\,. See also \cite{Adams2003}, pp. 252 Sect. 7.62\,.

\medskip

\noindent\textbf{Equivalence. }It is not so difficult to show that there exists an explicit constant $C_{s,d}>0$ such that
\begin{equation}\label{equiv.WH.Rd}
[f]_{W^{s,2}(\RR^d)}= C_{s,d}\left\|(-\Delta_{\RR^d})^{\frac{s}{2}}f\right\|_{\LL^2(\RR^d)}\,
=C_{s,d}\left\||\xi|^s \hat{f}\right\|_{\LL^2(\RR^d)}
\end{equation}
see e.g. Prop. 3.6 of \cite{NPV}\,, (in \cite{LM}, pg 59, Rem. 10.5 same claim, no explicit proof). Therefore, we can identify the spaces $H^s(\RR^d)=W^{s,2}(\RR^d)$ and equip them with the equivalent norms $\|\cdot\|_{H^s(\RR^d)}^2 \asymp \|\cdot\|_{W^{s,2}(\RR^d)}^2$\,, namely
\begin{equation}\label{Hs.equiv.norms}
\|\cdot\|_{H^s(\RR^d)}^2 \asymp \|\cdot\|_{W^{s,2}(\RR^d)}^2
\asymp \left\|f\right\|_{\LL^2(\RR^d)}^2+\left\|(-\Delta_{\RR^d})^{\frac{s}{2}}f\right\|_{\LL^2(\RR^d)}^2\,.
\end{equation}

\medskip

\noindent\textbf{Dual Spaces: $H^{-s}(\RR^d)$. }We can define the space $H^{-s}(\RR^d)$ as the dual space of $H^s(\RR^d)$\,, and, in view of the above discussion, we have that
\begin{equation}\label{H-s.Equivalence.Space}
H^{-s}(\RR^d)=\left(H^{s}(\RR^d)\right)^*=\left(W^{s,2}(\RR^d)\right)^*\,.
\end{equation}
We can characterize the space $H^{-s}(\RR^d)$ via the Fourier transform (even if $H^{-s}(\RR^d)$ is not a subset of $\LL^2(\RR^d)$\,, therefore we have to use the Fourier transform in $\mathcal{S}'(\RR^d)$\,, the space of tempered distributions):
\begin{equation}
H^{-s}(\RR^d)=\left\{f\in\mathcal{S}'(\RR^d)\,\big|\, \left(1+|\xi|^2\right)^{-\frac{s}{2}}\hat{f}\in \LL^2(\RR^d)\right\}\,,
\end{equation}
which is a Hilbert space with norm:
\begin{equation}\label{H-s.Equivalence.Norms}
\|f\|_{H^{-s}(\RR^d)}= \left\|\left(1+|\xi|^{2}\right)^{-\frac{s}{2}}\hat{f}\right\|_{\LL^2(\RR^d)}^2
\end{equation}
Notice that
\[
\|f\|_{H^{-s}(\RR^d)}= \left\|\left(1+|\xi|^{2}\right)^{-\frac{s}{2}}\hat{f}\right\|_{\LL^2(\RR^d)}^2
\le  \left\||\xi|^{-s}\hat{f}\right\|_{\LL^2(\RR^d)}^2 = \|(-\Delta)^{-\frac{s}{2}}f\|_{\LL^2(\RR^d)}\,,
\]
but the converse inequality is not true, therefore the norm on $H^{-s}$ is not equivalent to $\|(-\Delta)^{-\frac{s}{2}}\cdot\|_{\LL^2(\RR^d)}$\,. This issue is related to the dual of the spaces $\dot{H}^s(\RR^d)$\,, which we will not discuss here.

\medskip

\noindent\textbf{Sobolev-type Inequalities. }Let us recall the well-known fractional Sobolev inequality: there exists a positive constant $\mathcal{S}$ such that
\begin{equation}\label{Sob.Rd}
\|f\|_{\LL^{2^*}(\RR^d)}\le \mathcal{S} \left\|(-\Delta_{\RR^d})^{\frac{s}{2}}f\right\|_{\LL^2(\RR^d)}
= \mathcal{S} C_{s,d} [f]_{W^{s,2}(\RR^d)}\,,\qquad\mbox{for all $f\in H^s(\RR^d)$ and $2^*=\frac{2d}{d-2s}$}
\end{equation}
For a proof see e.g. \cite{LiebAnn} and  \cite{CotTav} or \cite{NPV}\,.

Notice that there is also a dual version of this Sobolev inequality, often called in literature the Hardy-Littlewood-Sobolev inequality:
\begin{equation}\label{HLS.Rd}
\|(-\Delta_{\RR^d})^{-s} g\|_{\LL^2(\RR^d)}\le \mathcal{S} \|g\|_{\LL^{q}(\RR^d)}\,,\qquad\mbox{for all $g\in H^{-s}(\RR^d)$ and $q=\frac{2d}{d+2s}$\,.}
\end{equation}
Here,
\[
(-\Delta_{\RR^d})^{-s} g(x):=c'_{d,s}\int_{\RR^d} \frac{g(y)}{|x-y|^{d-2s}}\dx
\]
A proof of \eqref{HLS.Rd} with optimal constant and explanation of the duality is given in \cite{LiebAnn}\,. We will explain more about this duality and the connection with the Green operator in Section \ref{Sect.Duality.Legendre}\,.

In the case of bounded domains, the situation is a bit more complicated, as we shall analyze in details in the rest of this section. From now on, we will adopt the notation $H^s(\RR^d)$ with the norm $\|\cdot\|_{H^s(\RR^d)}$ to indicate the Fractional Sobolev space on $\RR^d$ (any of the equivalent characterizations)\,, equipped with any of the equivalent norms in \eqref{Hs.equiv.norms}\,.

\subsection{The space $W^{s,p}(\Omega)$}

In the literature, fractional Sobolev-type spaces are also called Aronszajn, Gagliardo or Slobodeckij spaces, by the name of
the ones who introduced them, almost simultaneously (see \cite{Arosj, Ga, Slobo}).

We refer also to the books of Lions and Magenes \cite{LM} or Adams and Fournier \cite{Adams2003}

\noindent These fractional Sobolev spaces are usually defined through the so-called Gagliardo seminorm:
\begin{equation}
[f]_{W^{s,p}(\Omega)}:=\left(\int_\Omega\int_\Omega \frac{|f(x)-f(y)|^p}{|x-y|^{d+ps}}\dx\dy\right)^\frac{1}{p}\,,
\end{equation}
then we define the space
\begin{equation}\label{Wps.Spaces}
W^{s,p}(\Omega)=\left\{f\in\LL^p(\Omega)\,\big|\,[f]_{W^{s,p}(\Omega)}<+\infty \right\}
\end{equation}
which is a Banach space (Hilbert space for $p=2$), with the norm
\begin{equation}
\|f\|_{W^{s,p}(\Omega)}= \left([f]_{W^{s,p}(\Omega)}+\|f\|_{\LL^p(\Omega)}\right)^{\frac{1}{p}}
\end{equation}
On the other hand, we can define the subspace $W_0^{s,p}(\Omega)$ as
\begin{equation}\label{W0ps.Spaces}
W_0^{s,p}(\Omega)=\overline{C_c^{\infty}(\Omega)}^{\|\cdot\|_{W^{s,p}(\Omega)}}
\end{equation}
From now on, we will only deal with the case $p=2$\,.

\noindent\textbf{Traces. }When $s>1/2$ and the domain is sufficiently smooth (for instance Lipschitz), then the functions of $W^{s,2}(\Omega)$ has a trace in $\LL^2(\partial\Omega)$\,, namely the trace operator ${\rm Tr}: W^{s,2}(\Omega)\to \LL^2(\partial\Omega)$ is a linear and continuous operator. In this case, it is easy to check that the kernel of the operator ${\rm Tr}$ is ${\rm Ker}({\rm Tr})={\rm Tr }^{-1}0=W^{s,2}_0(\Omega)$\,, see \cite{Adams2003, LM, Mazya}. Notice that the spaces $W^{s,2}(\Omega)$ have no trace when $0<s\le 1/2$\,, therefore we can identify $W^{s,2}(\Omega)=W^{s,2}_0(\Omega)$ when $s\le 1/2$\,, as we shall see later.

\noindent\textbf{Extensions from $\Omega$ to $\RR^d$. }When the domain is regular enough (for instance Lipschitz)\,, then it has the so-called extension property, for example, see Theorem 5.4 of \cite{NPV}\,:

\begin{thm}Under the running assumptions, for any $s\in(0,1]$ and $p>1$\,, there exists an extension $\tilde{u}$ of $u\in W^{s,2}(\Omega)$ such that $\tilde{u}\in W^{s,2}(\RR^d)$\,, such that $\tilde{u}_{|\Omega}=u$ and
\begin{equation}\label{ext.Wps.Omega.Rd}
\|\tilde{u}\|_{W^{s,p}(\RR^d)}\le C \|u\|_{W^{s,p}(\Omega)}
\end{equation}
\end{thm}
Notice that the above inequality implies that the norm:
\begin{equation}\label{Wps.Omega.norm.2}
\|u\|_{W^{s,p}(\Omega),e}:=\inf_{\begin{subarray}{c} F\in W^{s,p}(\RR^d) \\ U=u \mbox{\footnotesize \;a.e. in $\Omega$}\end{subarray}}\normalsize\|U\|_{W^{s,p}(\RR^d)}
\end{equation}
is equivalent to $\|u\|_{W^{s,p}(\Omega)}$\,, since it is always true that $\|u\|_{W^{s,p}(\Omega)}\le \|\tilde{u}\|_{W^{s,p}(\RR^d)}$\,.
Therefore we can identify the functions in $W^{s,p}(\Omega)$ as the restriction to $\Omega$ of functions of $W^{s,p}(\RR^d)$\,, namely
\begin{equation}\label{Wsp.Omega.Caract.1}
W^{s,p}(\Omega)=\left\{f_{|\Omega}\;\big|\; f\in W^{s,p}(\RR^d)\right\}
\end{equation}
and we can endow the space with the norm $\|u\|_{W^{s,p}(\Omega),e}$ defined in \eqref{Wps.Omega.norm.2}\,.

Let us call $r_\Omega: W^{s,p}(\RR^d)\to W^{s,p}(\Omega)$ the restriction. By the above discussion, it is clear that $r_\Omega$ is surjective. Therefore,
\[
{\rm Ker}(r_\Omega)=\left\{f\in W^{s,p}(\RR^d)\;\big|\; r_\Omega(f)=f_{|\Omega}=0\right\}\,,
\]
therefore $r_\Omega$ gives the isomorphism
\begin{equation}\label{Wsp.quotient}
W^{s,p}(\Omega)= \quotient{W^{s,p}(\RR^d)}{{\rm Ker}(r_\Omega)}
\end{equation}

\noindent\textbf{Sobolev-type inequalities. }When the domain has the extension property, then $W^{2,s}(\Omega)$ inherits the Sobolev inequality \eqref{Sob.Rd} from $\RR^d$\,, via the extension, more precisely:
\begin{equation}\label{Sob.Omega}
\|f\|_{\LL^{q}(\Omega)}\le \mathcal{S} \|f\|_{W^{s,2}(\Omega)}\,,\qquad\mbox{for all $f\in W^{s,2}(\Omega)$ and $0<q\le 2^*=\frac{2d}{d-2s}$}
\end{equation}
See e.g. \cite{Mazya}\,. Moreover, the imbedding $W^{s,2}(\Omega)\hookrightarrow \LL^q(\Omega)$ is compact for all $1\le q< 2^*$\,. See also Thm 5.4 and 6.7 of \cite{NPV}\,.

\subsection{The spaces $H^s(\Omega)$, $H^s_0(\Omega)$ and $H^{\frac{1}{2}}_{00}$}\label{sec.7.3}
Following the book \cite{LM}\,, we can define the spaces $H^s(\Omega)$ by interpolation between $H^1(\Omega)$ and $\LL^2(\Omega)$\,, in analogy to the case of $\RR^d$\,:
\begin{equation}\label{Hs.Omega.Def.Interp}
H^s(\Omega)=\big[H^1(\Omega)\,,\,\LL^2(\Omega)\big]_{1-s}\,.
\end{equation}
According to this definition, this space is an Hilbert space with the natural norm given by the interpolation Theorem \ref{Thm.Jcont.Appendix}. We recall in the Appendix \ref{Appendix.Interp} some basic definitions and facts about interpolation of Banach spaces.

\noindent\textbf{Characterization by restriction/extension. }When the domain $\Omega$ is sufficiently regular (we need the extension property to hold), it is possible to characterize the space $H^s(\Omega)$ as the restriction to $\Omega$ of  functions of $H^s(\RR^d)$\,, more precisely:
\begin{equation}\label{Hs.Omega.Caract.1}
H^s(\Omega)=\left\{f_{|\Omega}\;\big|\; f\in H^{s}(\RR^d)\right\}
\end{equation}
see \cite{LM} Thm. 9.1\,. Moreover, we can take the following equivalent norm: (see \cite{LM} Thm. 9.2)
\begin{equation}\label{Hs.Omega.norm.1}
\|f\|_{H^s(\Omega),(1)}=\inf_{\begin{subarray}{c} F\in H^s(\RR^d) \\ F=f \mbox{\footnotesize \;a.e. in $\Omega$}\end{subarray}}\normalsize\|F\|_{H^s(\RR^d)}\,.
\end{equation}
notice that $\|\cdot\|_{H^s(\RR^d)}$ denotes any equivalent norm on $H^s(\RR^d)$\,, as in \eqref{Hs.equiv.norms}\,. We can rephrase the above characterization in terms of the extensions from $\Omega$ to $\RR^d$, as we have done above for the case $W^{s,p}$\,.

We can now define the space $H^s_0(\Omega)$ by
\begin{equation}\label{def.Hs0.Omega.1}
H_0^s(\Omega)=\overline{C_c^{\infty}(\Omega)}^{\|\cdot\|_{H^s(\Omega)}}
\end{equation}
as it has been done in \cite{LM} pg. 60, sect 11.1. Excluding exceptional values of $s\in (0,1)$ as explained below, there another equivalent definition given via interpolation, namely Theorem 11.6 of \cite{LM}, which states that, (notice that here $s_1$ and $s_2$ are real numbers that can be greater than one)
\begin{equation}\label{def.Hs0.Omega.2}
H^s_0(\Omega):= \big[H^{s_1}_0(\Omega)\,,\,H^{s_2}_0(\Omega)\big]_\theta\,,\;\mbox{for any $0\le s_1\,,\, s_2 $ such that $s=(1-\theta)s_1+\theta s_2\neq \mbox{integer}+\frac{1}{2}$\,.}
\end{equation}
According to this definition, this space is an Hilbert space with the norm given by the interpolation Theorem \ref{Thm.Jcont.Appendix}. We recall in the Appendix \ref{Appendix.Interp} some basic definitions and facts about interpolation of Banach spaces.

\subsection{Properties of the spaces $H^s(\Omega)$ and $H^s_0(\Omega)$. }\label{sect.prop.Hs} We address here the properties of two related spaces that appear often in the analysis.

\noindent \textbf{Density and difference }
To this end, we recall the result by Lions and Magenes.

\begin{thm}[Theorem 11.1 of \cite{LM}]
If the bounded domain $\Omega$ is smooth enough, then $C_c^{\infty}(\Omega)$ is dense in $H^s(\Omega)$ if and only if $0< s\le 1/2$; in this case we have that $H^s_0(\Omega)= H^s(\Omega)$.  If $s>1/2$\,, then $H^s_0(\Omega)\subset H^s(\Omega)$ and the inclusion is strict.
\end{thm}

\noindent \textbf{Traces. }Theorem 9.4 of \cite{LM} proves that the trace operator is well defined from ${\rm Tr}:H^s(\Omega)\to \LL^2(\partial\Omega)$ only when $s>1/2$\,. In this case, it is easy to check that the kernel of the operator is ${\rm Ker}({\rm Tr})={\rm Tr }^{-1}0=H^s_0(\Omega)$\,. Therefore, one can characterize the functions of $H^s_0(\Omega)$ when $1/2<s\le 1$ by:
\begin{equation}\label{Hs0.charact.2.trace}
H^s_0(\Omega)=\left\{f\in H^s(\Omega)\;\big|\; u=0\quad\mbox{in $\LL^2(\partial\Omega)$}\right\}\,.
\end{equation}

\noindent \textbf{Weak traces. Hardy-type inequalities. }On the other hand, when $0\le s< 1/2$\,, $H^s(\Omega)=H^s_0(\Omega)$ and the trace is not well defined, but there is at least a weaker concept, which basically states that a function of $H^s(\Omega)$ divided by the $s$-power of the distance to the boundary is  in $\LL^2(\Omega)$. Indeed, let $\phi$ be a $C^\infty$ extension of the distance to the boundary, namely a $\phi\in C^\infty(\overline{\Omega})$\,, with $\phi>0$ on $\Omega$ and $\phi=0$ on $\partial\Omega$\,, and moreover $\phi(x_0)/\dist(x_0,\partial\Omega)=\mbox{const.}$ for all $x_0\in\partial\Omega$\,. Then the operator ``multiplication by $\phi^{-s}$\,'' is linear and continuous from
\[
\begin{array}{ll}
\phi^{-s}:  H^s(\Omega)\to \LL^2(\Omega)& \qquad\mbox{when $0<s<\frac{1}{2}$}\\[2mm]
\phi^{-s}:  H_0^s(\Omega)\to \LL^2(\Omega) &\qquad\mbox{when $\frac{1}{2}<s\le 1$}\\
\end{array}
\]
see Thm 11.2 and Thm 11.3 of \cite{LM} respectively. Namely, there exists a constant $C>0$ such that for all $f\in H_0^s(\Omega)$: (recall that $H^s=H_0^s$ when $s\le 1/2$)
\begin{equation}\label{hardy.like}
\left\|\frac{f}{\phi^s}\right\|_{\LL^2(\Omega)}\le C\|f\|_{H^s(\Omega)}\,.
\end{equation}
The above inequality is indeed an Hardy-type inequality, cf. for example \cite{Davies1,Davies2}\,.
The case $s=1/2$ is a special case, as we shall discuss later.

\noindent\textbf{The extension by zero outside $\Omega$. }When the domain $\Omega$ is sufficiently regular, we already know that there is an extension of functions of $H^s(\Omega)$ to $H^s(\RR^d)$\,, cf. \eqref{Hs.Omega.Caract.1}\,. We are now interested to know when such extension is the trivial one, namely the extension by zero outside $\Omega$\,. The answer to this question is given by Thm. 11.4 of \cite{LM}: denote by $E_0(f)$ the zero extension of $f$ on $\RR^d\setminus\Omega$\,,then we have that the linear operator $E_0$ is continuous from
\[
\begin{array}{lll}
E_0:  H^s(\Omega)\to H^s(\RR^d)& \qquad\mbox{when $0<s<\frac{1}{2}$} \\[2mm]
E_0:  H_0^s(\Omega)\to H^s(\RR^d) &\qquad\mbox{when $\frac{1}{2}<s\le 1$}\\
\end{array}
\]
namely there exists a constant $C>0$ such that for all $f\in H_0^s(\Omega)$: (recall that $H^s(\Omega)=H_0^s(\Omega)$ when $s\le 1/2$)
\[
\|E_0(f)\|_{H^s(\RR^d)}\le C \|f\|_{H^s(\Omega)}\,.
\]
Moreover, we have now a more explicit expression of the norm $\|f\|_{H^s(\Omega),(1)}$ defined in \eqref{Hs.Omega.norm.1}\,, with the help of the zero-extension, for any $f\in H_0^s(\Omega)$:
\begin{equation}\label{Hs.Omega.norm.2}
\|f\|_{H^s(\Omega),e}=\inf_{\begin{subarray}{c} F\in H^s(\RR^d) \\ F=f \mbox{\footnotesize \;a.e. in $\Omega$}\end{subarray}}\normalsize\|F\|_{H^s(\RR^d)}= \|E_0(f)\|_{H^s(\RR^d)}
\end{equation}
notice that $\|\cdot\|_{H^s(\RR^d)}$ denotes any equivalent norm on $H^s(\RR^d)$\,, as in \eqref{Hs.equiv.norms}\,.\\
The case $s=1/2$ is a special case, as we shall discuss later.

Let us call $r_\Omega: H^s(\RR^d)\to H^s(\Omega)$ the restriction. By the above discussion, it is clear that $r_\Omega$ is surjective. We have
\[
{\rm Ker}(r_\Omega)=\left\{f\in H^s(\RR^d)\;\big|\; r_\Omega(f)=f_{|\Omega}=0\right\}\,,
\]
therefore $r_\Omega$ gives the isomorphism
\begin{equation}\label{Hs.quotient}
H^s(\Omega)= \quotient{H^s(\RR^d)}{{\rm Ker}(r_\Omega)}
\end{equation}

\subsection{Isomorphism between $W^{s,2}(\Omega)$ and $H^s(\Omega)$ spaces} We shall see now that the different definitions of $H^s$and $W^{s,2}$ spaces are the same, through the extension to the whole space. Of course the domain $\Omega$ need to be sufficiently regular, namely the extension property shall holds true. More precisely, let us recall the isomorphism  \eqref{Wsp.quotient} and recall that $r_\Omega$ is the restriction, so that
\begin{equation}\label{equiv.WH.Omega}
W^{s,2}(\Omega)= \quotient{W^{s,2}(\RR^d)}{{\rm Ker}(r_\Omega)}= \quotient{H^s(\RR^d)}{{\rm Ker}(r_\Omega)}=H^s(\Omega)
\end{equation}
the isomorphism in the middle holds since in $\RR^d$ we already know that $W^{s,2}(\RR^d)=H^s(\RR^d)$ as consequence of \eqref{equiv.WH.Rd} and \eqref{Hs.equiv.norms}\,.

Moreover, by transposition we get that $r_\Omega^*$ gives the isomorphism of duals
\begin{equation}\label{Hs.quotient.dual.1}\begin{split}
\left(H^s(\Omega)\right)^*
&= \left(\quotient{H^s(\RR^d)}{{\rm Ker}(r_\Omega)}\right)^*\\
&=\left\{f\in H^{-s}(\RR^d)  \;\big|\; \int_{\RR^d}f\psi\dx=0\qquad\mbox{for all $\psi\in C_c^\infty(\RR^d)$ with $\psi=0$ on $\RR^d\setminus\Omega$} \right\}
\end{split}
\end{equation}
Defining then
\begin{equation}\label{Hs.quotient.dual.2}
H^{-s}_{\overline{\Omega}}(\RR^d)= \left\{f\in  H^{-s}(\RR^d) \;\big|\; \supp{f}\subseteq\overline{\Omega} \right\}
\end{equation}
we have that $r_\Omega^*$ gives the isomorphism between $\left(H^s(\Omega)\right)^* $ and $H^{-s}_{\overline{\Omega}}(\RR^d)$\,.

\noindent The above discussion is similar to Remarks 12.4 and 12.5 of \cite{LM}\,.   Summing up, we obtain the following
\begin{thm}\label{thm.Hs-Ws}
Let $\Omega$ be a bounded regular domain of $\RR^d$, then $H^s(\Omega)=W^{s,2}(\Omega)$ and the same holds for the dual spaces.
\end{thm}

\subsection{The space $H^{\frac{1}{2}}_{00}(\Omega)$. }The case $s=1/2$ is special, as one can deduce from the above discussion. The problem is how to define, if possible, a weak concept of trace in this space. This can be done by interpolation, indeed we define
\begin{equation}\label{def.H1200.1}
H^{\frac{1}{2}}_{00}(\Omega):= \big[H^{s_1}_0(\Omega)\,,\,H^{s_2}_0(\Omega)\big]_\theta\,,\qquad\mbox{for any $ 0\le s_2<\frac{1}{2}<s_1$\,, such that $(1-\theta)s_1+\theta s_2=\frac{1}{2}$\,.}
\end{equation}
According to this definition, this space is an Hilbert space with the norm given by the interpolation Theorem \ref{Thm.Jcont.Appendix}. We recall in the Appendix \ref{Appendix.Interp} some basic definitions and facts about interpolation of Banach spaces.

One may wonder whether the above definition depends on the particular choice of $s_1$ or $s_2$: Theorem 11.7 of \cite{LM}\,, gives a negative answer to this question: the space $\big[H^{s_1}_0(\Omega)\,,\,H^{s_2}_0(\Omega)\big]_\theta$ is independent on $s_i$ satisfying the relation \eqref{def.H1200.1}\,. Moreover, the following characterization holds:
\begin{equation}\label{def.H1200.2}
H^{\frac{1}{2}}_{00}(\Omega):= \left\{f\in H^{\frac{1}{2}}_0(\Omega)\;\big|\; \phi^{-\frac{1}{2}}f\in \LL^2(\Omega)\right\}\,,
\end{equation}
where $\phi$ is a $C^\infty$ extension of the distance to the boundary, as in \eqref{hardy.like}\,. Moreover, the interpolation norm is equivalent to
\begin{equation}\label{norm.H1200.1}
\|f\|_{H^{\frac{1}{2}}_{00}(\Omega)}:= \left(\|f\|^2_{H^{\frac{1}{2}}(\Omega)}
    +\left\|\phi^{-\frac{1}{2}}f\right\|_{\LL^2(\Omega)}^2\right)^{\frac{1}{2}}\,,
\end{equation}
Finally, the space $H^s_{00}(\Omega)$ is strictly included in $H^s_{0}(\Omega)$ with a strictly finer topology.
Roughly speaking, we are considering the closed subspace of $H^{\frac{1}{2}}_0(\Omega)$ such that the Hardy inequality holds true.

\subsection{Relation with the dual spaces $W^{-s,2}$ and $H^{-s}(\Omega)$}
Let us define the dual spaces, for all $0<s\le 1$
\begin{equation}\label{WH.duals}
H^{-s}(\Omega)=\left(H^s_0(\Omega)\right)^*\qquad\mbox{and}\qquad W^{-s,2}(\Omega)=\left(W^{s,2}_0(\Omega)\right)^*\,.
\end{equation}
We have seen in the previous section the isomorphism between spaces $W^{s,2}(\Omega)$ and $H^s(\Omega)$\,, moreover we know that when $0<s<1/2$ we have $H^s_0(\Omega)=H^s(\Omega)$\,, so that by \eqref{Hs.quotient.dual.2} we obtain
\begin{equation}\label{WH.duals<12}
H^{-s}(\Omega)=W^{-s,2}(\Omega)=H^{-s}_{\overline{\Omega}}(\RR^d)=\left\{f\in H^{-s}(\RR^d)\;\big|\; \supp{f}\subseteq\overline{\Omega}  \right\} \qquad\mbox{for all $0<s<\frac{1}{2}$}\,.
\end{equation}
We can equip $H^{-s}(\Omega)$ with the following norm:
\begin{equation}\label{H-s.Omega.norm.1}
\|f\|_{H^{-s}(\Omega)}=\|F\|_{H^{-s}(\RR^d)}=\|E(f)\|_{H^{-s}(\RR^d)}
\end{equation}
where $E$ is the inverse of the isomorphism $r_\Omega^*$\,.

On the other hand, when $1/2<s\le 1$ the spaces $H_0^s$ are strictly included in $H^s$\,, therefore the dual $H^{-s}(\Omega)$ is a projection of $\left(H^s(\Omega)\right)^*=H^{-s}_{\overline{\Omega}}(\RR^d)$\,. When $s=1/2$ one has to notice that the dual spaces of $H_{00}^{\frac{1}{2}}(\Omega)$is different from the dual of $H^{\frac{1}{2}}(\Omega)=H_0^{\frac{1}{2}}(\Omega)$.

In view of the characterization \eqref{def.H1200.2} of  $H_{00}^{\frac{1}{2}}(\Omega)$, we can characterize its dual as follows
\begin{equation}\label{def.H1200.dual}
\left(H^{\frac{1}{2}}_{00}(\Omega)\right)^*:= \left\{f\in\mathcal{D}'(\Omega) \;\big|\; f=f_0+f_1\quad\mbox{with}\quad f_0\in H^{-\frac{1}{2}}(\Omega)\quad\mbox{and}\quad \phi^{\frac{1}{2}}f_1\in \LL^2(\Omega)\right\}\,,
\end{equation}
We shall finally remark that the function $\phi^{\varepsilon-1}\in \left(H^{\frac{1}{2}}_{00}(\Omega)\right)^*$ for all $\varepsilon>0$\,, where $\phi$ is a $C^\infty$ extension of the distance to the boundary, as in \eqref{hardy.like}\,. For more details, see Remarks 12.1 of \cite{LM}\,.

%%%%%%%%%%%%%%%%%%%%%%%%%%%%%%%%%%%%%%%%%%%%%%%%%%%%%%%%%%%%%%%%%%%%%%

\subsection{A dual approach to functional inequalities via Legendre transform}\label{Sect.Duality.Legendre} It would be interesting to prove the validity of suitable functional inequalities of Sobolev type for the norm of $H$ or $H^*$\,.   We begin by proving an equivalence result.
\begin{prop}\label{prop.equiv.ineq}The following inequalities are equivalent:\\
$(I)$ For all $f\in H$ and all $0<q\le 2^*=\frac{2d}{d-2s}$
\begin{equation}\label{Sob.H}
\|f\|_{\LL^q(\Omega)}\le C \|\AM f\|_{\LL^2(\Omega)}=\|f\|_H\,,
\end{equation}
$(II)$ For all $g\in H^*$ and all $q'\ge (2^*)'=\frac{2d}{d+2s}$
\begin{equation}\label{Sob.H*}
\|g\|_{H^*}=\|\AIM g\|_{\LL^2(\Omega)}\le C \|g\|_{\LL^{q'}(\Omega)}\,.
\end{equation}
\end{prop}
\noindent {\sl Proof.~}The fact that the above two inequalities are equivalent (and that the constant $C>0$ is the same) can be seen trough the Legendre transform: let $N:\,H \to \RR$ be a convex functional on $H$\,, define the Legendre transform of $N$ as the functional $N^*:H^*\to\RR$ such that
\[
N^*[g]:=\sup_{f\in D}\big(\langle f\,,\,g\rangle_{H,H^*} - N[f]\big)
\]
where $D\subset H$ can be a dense subset of $H$\,.

\noindent The statement of the proposition is equivalent to the following claim.  \textsl{For all $f\in H\cap \LL^q(\Omega)$ and $g\in H^*\cap \LL^{q'}(\Omega)$ with $\dfrac{1}{q}+\dfrac{1}{q'}=1$\,, we have
\[
N_1[f]=\frac{1}{2}\|f\|_{\LL^q(\Omega)}^2\qquad\mbox{if and only if}\qquad N_1^*[g]=\frac{1}{2}\|g\|^2_{\LL^{q'}(\Omega)}
\]
and (proof by series expansions)}
\[
N_2[f]=\frac{C}{2}\|\AM f\|_{\LL^2(\Omega)}^2\qquad\mbox{if and only if}\qquad N_2^*[g]=\frac{2}{C}\|\AIM g\|_{\LL^2(\Omega)}^2
\]
Indeed, let us compute
\[\begin{split}
N_1^*[g]
&=\sup_{f\in H\cap \LL^q(\Omega)}\left(\int_\Omega f\,g\dx-\frac{1}{2}\left(\int_\Omega |f|^q\dx\right)^{\frac{2}{q}}\right)\\
    &=\sup_{t\ge 0}\sup_{\begin{subarray}{c} f\in H\cap \LL^q(\Omega)\\ \|f\|_{\LL^{q}(\Omega)}=t\end{subarray}}
        \left(\int_\Omega f\,g\dx-\frac{1}{2}\left(\int_\Omega |f|^q\dx\right)^{\frac{2}{q}}\right)
 =\sup_{t\ge 0}\left(t\|g\|_{\LL^{q'}(\Omega)}-\frac{t^2}{2}\right)=\frac{1}{2}\|g\|_{\LL^{q'}(\Omega)}^2\,.
\end{split}
\]
On the other hand,
\[\begin{split}
N_2^*[g]
&=\sup_{f\in H}\left(\int_\Omega f\,g\dx-\frac{C}{2}\|f\|^2_{H}\right)\\
    &=\sup_{t\ge 0}\sup_{\begin{subarray}{c} f\in H\\ \|f\|_{H}=t\end{subarray}}
        \left(\int_\Omega f\,g\dx-\frac{C}{2}\|f\|_H^2\right)
 =\sup_{t\ge 0}\left(t\|g\|_{H^*}-\frac{C}{2}t^2\right)=\frac{1}{2 C}\|g\|_{H^*}^2\mbox{\,.\qed}
\end{split}
\]
Legendre's duality then implies that inequality \eqref{Sob.H} holds if and only if inequality \eqref{Sob.H*} holds: namely
\[
N_1[f]\le N_2[f]\qquad\mbox{for all $f\in H$\,,}\qquad\mbox{if and only if}\qquad N_2^*[g]\le N_1^*[g]\qquad\mbox{for all $g\in H^*\,.$}
\]
recall that a priori the above inequalities hold in $H\cap\LL^q(\Omega)$) and in $H^*\cap\LL^{q'}(\Omega)$ respectively. Then we shall recall that both subspaces are dense in $H$ and $H^*$ respectively, therefore the inequalities can be extended to $H$ and $H^*$  respectively, by density.\qed

\subsubsection{Sobolev type inequalities}
We now show that suitable Sobolev inequalities holds true as a consequence of the bounds on the Green function \eqref{Gree.est.1}, namely
$ 0\le G_{\Omega}(x,x_0)\le c_{1,\Omega}\,|x-x_0|^{-(d-2s)}$, for all $x,x_0 \in \Omega$\,.
\begin{thm}The following inequality holds for all $f\in H$ and all $0<q \le 2^*=\frac{2d}{d-2s}$,
\begin{equation}\label{Sob.H.1}
\|f\|_{\LL^q(\Omega)}\le C \|\AM f\|_{\LL^2(\Omega)}=\|f\|_H\,,
\end{equation}
Moreover, it is equivalent to the dual inequality: for all $g\in H^*$ and all $q'\ge (2^*)'=\frac{2d}{d+2s}$,
\begin{equation}\label{Sob.H*.1}
\|g\|_{H^*}=\|\AIM g\|_{\LL^2(\Omega)}\le C \|g\|_{\LL^{q'}(\Omega)}\,.
\end{equation}
\end{thm}
\noindent {\sl Proof.~}By Proposition \ref{prop.equiv.ineq} we get that inequality \ref{Sob.H.1} holds if and only if inequality \eqref{Sob.H*.1} holds.   We are going to prove the dual inequality \eqref{Sob.H*.1} using the Green function estimates \eqref{Gree.est.1}, and using the corresponding Hardy-Littlewood-Sobolev inequality on the whole space $\RR^d$\,, which reads:
\begin{equation}\label{Sob.H*.2}
\|G_{\RR^d} g\|_{\LL^2(\RR^d)}\le \mathcal{S} \|g\|_{\LL^{q'}(\RR^d)}\,,\qquad\mbox{for $q'=(2^*)'=\frac{2d}{d+2s}$\,.}
\end{equation}
where we recall that the Green function of the fractional Laplacian on $\RR^d$ takes the explicit form $G^s_{\RR^d}(x,y)=c_d\,/\,|x-y|^{d-2s}$\,, therefore thanks to the upper estimates \eqref{Gree.est.1} on the Green function we have for all $x,x_0 \in \Omega$:
\begin{equation}\label{Gree.est.2}
0\le G_{\Omega}(x,x_0)\le \frac{c_{1,\Omega}}{|x-x_0|^{d-2s}}
\le C_1\,G_{\RR^d}(x,x_0)
\end{equation}
and we shall extend this inequality to zero for $x,x_0\in \RR^d\setminus \Omega$\,.\\
Finally, consider $g\in H^*$\,, and call $g_0$ its zero extension in $\RR^d\setminus \Omega$\,, then one has:
\[\begin{split}
\|g\|_{H^*}^2&=\int_\Omega\left(\int_\Omega G_{\Omega}(x,y)g(y) \dy\right)^2\dx
    \le_{(a)}  C_1^2\,\int_\Omega\left(\int_\Omega G_{\RR^d}(x,y) g(y) \dy\right)^2\dx\\
&=_{(b)} C_1^2\, \left\|\widetilde{\AI_{\RR^d} g_0}\right\|_{\LL^2(\RR^d)}^2
\le_{(c)} C_1^2\, \left\|\AI_{\RR^d} g_0\right\|_{\LL^2(\RR^d)}^2
\le_{(d)} C_1^2\,\mathcal{S}^2\|g_0\|_{\LL^{q'}(\RR^d)}^2
= C_2^2\|g\|_{\LL^{q'}(\Omega)}^2
\end{split}
\]
where in $(a)$ we have used inequality \eqref{Gree.est.2}\,, in $(b)$ we have defined $\widetilde{\AIM_{\RR^d} g_0}$ as the zero truncation outside $\Omega$ of the function $x\mapsto \int_\Omega G_{\RR^d}(x,y) g(y) \dy$\,, notice that \textit{this is not the solution $F$ of the equation $(-\Delta_{\RR^d})^s F = g_0$ in the whole space $\RR^d$ \,, corresponding to $g_0$\,, which is the zero extension outside $\Omega$ of $g\in H$\,. }Indeed, since $g(y)\, G_{\RR^d}(x,y)\ge 0$\,, we have
\[
\int_\Omega G_{\RR^d}(x,y) g(y) \dy \le \int_{\RR^d} G_{\RR^d}(x,y) g_0(y) \dy =F(x)
\]
Anyway the above inequality, shows that inequality $(c)$ is true, provided $F\in \LL^2(\RR^d)$\,. Indeed, since $\supp(g_0)=\overline{\Omega}$, then $F$ is not necessarily in $\LL^1(\RR^d)$\,, but it can be shown to belong to $\LL^2(\RR^d)$ whenever $g_0\in \LL^{q'}(\RR^d)$ with $q'=(2^*)'=\frac{2d}{d+2s}$\,, by the Hardy-Littlewood-Sobolev inequality \eqref{Sob.H*.2}\,, namely
\[
\|F\|_{\LL^2(\RR^d)}=\|\AIM_{\RR^d} g_0\|_{\LL^2(\RR^d)}\le \mathcal{S} \|g_0\|_{\LL^{q'}(\RR^d)}\,.
\]
Therefore, inequality $(c)$ and $(d)$ hold true. Finally, we have proved that the dual Sobolev inequality \eqref{Sob.H*.1} holds, therefore also the Sobolev inequality \eqref{Sob.H.1} holds.\qed

\subsection{Some technical facts}\label{Proofs.Sect.211}

\noindent\textbf{Proofs of section \ref{sect.H.spaces}}

\noindent {\sl Proof.~}First we prove that $\|F\|_{H^*}\le\left(\sum_{k=1}^{\infty}\lambda_k^{-1}\hat{F}_k^2\right)^{\frac{1}{2}}$\,, indeed
\begin{equation}\label{ineq.le.H*}\begin{split}
\|F\|_{H^*}
&=\sup_{\begin{subarray}{c}  g\in H \\ \|g\|_{H}\le 1\end{subarray}}\sum_{k=1}^{\infty}\hat{F}_k\,\hat{g}_k
=\sup_{\begin{subarray}{c}  g\in H \\ \|g\|_{H}\le 1\end{subarray}}
    \sum_{k=1}^{\infty}\lambda_k^{-\frac{1}{2}}\hat{F}_k\,\lambda_k^{\frac{1}{2}}\hat{g}_k\\
&\le\sup_{\begin{subarray}{c}  g\in H \\ \|g\|_{H}\le 1\end{subarray}}
    \left(\sum_{k=1}^{\infty}\lambda_k^{-1}\hat{F}_k^2\right)^{\frac{1}{2}}
    \left(\sum_{k=1}^{\infty}\lambda_k\hat{g}_k^2\right)^{\frac{1}{2}}
\le\left(\sum_{k=1}^{\infty}\lambda_k^{-1}\hat{F}_k^2\right)^{\frac{1}{2}}
    \sup_{\begin{subarray}{c}  g\in H \\ \|g\|_{H}\le 1\end{subarray}}\|g\|_{H}
=\left(\sum_{k=1}^{\infty}\lambda_k^{-1}\hat{F}_k^2\right)^{\frac{1}{2}}
\end{split}
\end{equation}
On the other hand, to prove that $\|F\|_{H^*}\ge \left(\sum_{k=1}^{\infty}\lambda_k^{-1}\hat{F}_k^2\right)^{\frac{1}{2}}$ we choose a special $g$ of the form
\[
\underline{g}(x)=\frac{\sum_{k=1}^{\infty}\lambda_k^{-1}\hat{F}_k\Phi_k(x)}{\left(\sum_{k=1}^{\infty}\lambda_k^{-1}\hat{F}_k^2\right)^{\frac{1}{2}}}\,.
\]
Then, recalling that $\Phi_k$ is an ortonormal basis of $\LL^2(\Omega)$\,, we get
\[
\hat{\underline{g}}_n=\int_\Omega \underline{g}(x)\Phi_n(x)\dx
=\frac{\sum_{k=1}^{\infty}\lambda_k^{-1}\hat{F}_k
\int_\Omega \Phi_k(x)\Phi_n(x)\dx}{\left(\sum_{k=1}^{\infty}\lambda_k^{-1}\hat{F}_k^2\right)^{\frac{1}{2}}}
=\frac{\lambda_n^{-1}\hat{F}_n}{\left(\sum_{k=1}^{\infty}\lambda_k^{-1}\hat{F}_k^2\right)^{\frac{1}{2}}}\,,
\]
so that $\|\underline{g}\|_H=1$\,, namely:
\[
\|\underline{g}\|_H=\left(\sum_{n=1}^{\infty}\lambda_n\hat{\underline{g}}_n^2\right)^{\frac{1}{2}}
=\left(\frac{\sum_{n=1}^{\infty}\lambda_n\,\lambda_n^{-2}\hat{F}_n^2}{\sum_{k=1}^{\infty}\lambda_k^{-1}\hat{F}_k^2}\right)^{\frac{1}{2}}=1
\]
and we have
\begin{equation}\label{ineq.ge.H*}\begin{split}
\|F\|_{H^*}
&=\sup_{\begin{subarray}{c}  g\in H \\ \|g\|_{H}\le 1\end{subarray}}\langle F, g \rangle_{H^*,H}
\ge  \langle F, \underline{g} \rangle_{H^*,H}
=\sum_{n=1}^{\infty}\hat{F}_n\,\hat{\underline{g}}_n
=\frac{\sum_{n=1}^{\infty}\lambda_n^{-1}\hat{F}_n^2}{\left(\sum_{k=1}^{\infty}\lambda_k^{-1}\hat{F}_k^2\right)^{\frac{1}{2}}}
=\left(\sum_{k=1}^{\infty}\lambda_k^{-1}\hat{F}_k^2\right)^{\frac{1}{2}}
\end{split}
\end{equation}
Combining \eqref{ineq.le.H*} and \eqref{ineq.ge.H*} we obtain \eqref{def.H*.norm.2}\,.\qed

\noindent\textbf{Proof of \ref{Int.by.Parts}. }
follows by the orthogonality property of the eigenfunctions $\Phi_k$, indeed
\[\begin{split}
\int_\Omega A^{1/2} f \,  A^{1/2} g\dx
&=\int_\Omega \sum_{k=1}^{\infty}\lambda_k^{\frac{1}{2}}\hat{f}_k\Phi_k(x)\,\sum_{n=1}^{\infty}\lambda_n^{\frac{1}{2}}\hat{g}_n\Phi_n(x)\dx\\
&=\sum_{k=1}^{\infty}\sum_{n=1}^{\infty}\lambda_n^{\frac{1}{2}}\lambda_k^{\frac{1}{2}}\hat{g}_n\hat{f}_k\int_\Omega \Phi_k(x)\,\Phi_n(x)\dx
=\sum_{k=1}^{\infty}\sum_{n=1}^{\infty}\lambda_n^{\frac{1}{2}}\lambda_k^{\frac{1}{2}}\hat{g}_n\hat{f}_k\delta_{k,n}
= \sum_{k=1}^{\infty}\lambda_k\hat{f}_k\hat{g}_k
\end{split}
\]
so that
\[\begin{split}
\int_\Omega f  A  g\dx
&= \sum_{k=1}^{\infty}\lambda_k\hat{g}_k\int_\Omega f(x) \,\Phi_k(x)\dx
    =\sum_{k=1}^{\infty}\lambda_k\hat{f}_k\hat{g}_k
    = \int_\Omega A^{1/2} f \,  A^{1/2} g\dx=\int_\Omega g  A  f\dx
\end{split}
\]
As a consequence of this formula, we can express the scalar product and the norm of $H$ in the following equivalent ways, for any $f,g\in H$\,:
\begin{equation}\label{product.H}
\langle f, g \rangle_H = \sum_{k=1}^{\infty}\lambda_k\hat{f}_k\,\hat{g}_k
=\int_\Omega A^{1/2} f \,  A^{1/2} g\dx
=\int_\Omega f  A  g\dx=\int_\Omega g  A  f\dx\,.
\end{equation}
We thus get three equivalent expressions for the norm:
\begin{equation}\label{norm.H}
\|f\|_H = \sum_{k=1}^{\infty}\lambda_k\hat{f}_k^2
=\left\|A^{1/2} f\right\|_{\LL^2(\Omega)}^2
=\int_\Omega f  A  f\dx\,.
\end{equation}

\

\section{Appendix 2. Reminder about interpolation spaces}\label{Appendix.Interp}

We recall the basic definitions and results about interpolation between Banach spaces, following the book \cite{Adams2003} Chapter 7. We will use the so-called J-method. Let $X_0,X_1$ be Banach spaces, we know that $X_0\cap X_1\subseteq X_0,X_1\subseteq X_0+X_1$\,, we recall that we have the following norms:
\[\begin{split}
\|u\|_{X_0\cap X_1}&=\max\{\|u\|_{X_0}\,,\,\|u\|_{X_1}\}\\
\|u\|_{X_0+X_1}&=\inf\{\|u_0\|_{X_0}+\|u_1\|_{X_1}\;\big|\; u=u_0+u_1\quad\mbox{with}\; u_0\in X_0\,,\,u_1\in X_1\}
\end{split}
\]

Define moreover the space (Bochner integral)
\[
\LL^q_*:=\LL^1\left((0,+\infty);\frac{dt}{t}\,:\,X_0+X_1\right)\qquad\mbox{for any $q\ge 1$\,.}
\]
\noindent\textbf{The J-method. }Define the $J(t,u)$-norm by
\begin{equation}\label{def.J}
J(t,u)=\max\{\|u\|_{X_0}\,,\,t\,\|u\|_{X_1}\}\,.
\end{equation}
The $J$-norm is clearly equivalent to $\|u\|_{X_0\cap X_1}$\,. If $0\le \theta\le 1$ and $1\le q\le \infty$ we denote by
\[
X_{\theta,q}:=\big[X_0,X_1\big]_{\theta,q}
\]
the space of all $u\in X_0+X_1$ such that
\[
u=\int_0^{+\infty}f(t)\frac{\dt}{t}
\]
for some $f\in \LL^1_*$ having values in $X_0\cap X_1$\,, and such that the function
\[
t\mapsto \frac{J(t,f(t))}{t^\theta}\in \LL^q_*
\]
In the case $q=2$\,, we will simplify the notation as follows:
\[
X_{\theta,2}:=\big[X_0,X_1\big]_{\theta,2}=\big[X_0,X_1\big]_{\theta}=X_\theta\,.
\]
Then we have the following Theorem:
\begin{thm}[The J-Method]\label{Thm.Jcont.Appendix}
If either $1<q\le \infty$ and $0<\theta<1$ or $q=1$ and $0\le\theta\le 1$\,, then $X_{\theta,q}=\big[X_0,X_1\big]_{\theta,q}$ is a non-trivial Banach space with the norm
\begin{equation}\label{norm.interp.cont.1}
\|u\|_{X_{\theta,q},J}=\inf_{f\in S(u)}\|f\|_{\LL^q_*}=\inf_{f\in S(u)}\left(\int_0^\infty\left(t^{-\theta}J(t,f(t))\right)^q\frac{\dt}{t}\right)^{\frac{1}{q}}\qquad\mbox{if $q<\infty$}\,,
\end{equation}
where
\begin{equation}
S(u)=\left\{f\in \LL^1_*\;\big|\; u=\int_0^{+\infty}f(t)\frac{\dt}{t}\right\}
\end{equation}
Furthermore,
\[
\|u\|_{X_0\cap X_1}\le \|u\|_{\theta,q,J}\le\|u\|_{X_0+X_1}
\]
so that $X_0\cap X_1\hookrightarrow \big[X_0,X_1\big]_{\theta,q}\hookrightarrow X_0+X_1$ with continuous injections, so that $\big[X_0,X_1\big]_{\theta,q}$ is an intermediate space between $X_0$ and $X_1$\,.
\end{thm}
\noindent {\sl Proof.~}This is  Thm. 7.13 of \cite{Adams2003}, proof at page 211.\qed

We are interested in a discrete version of the above theorem, but in a slightly more general form than the one given in Thm 7.14 of \cite{Adams2003}\,.

\begin{thm}[The Discrete version of the J-Method]\label{Thm.Jdiscr.Appendix}
Let $\mu_k$ be an increasing sequence $0<\mu_k<\mu_{k+1}\to +\infty$\,, such that $0<\mu_{k+1}/\mu_k\le \Lambda_0<\infty$\,. Let either $1<q\le \infty$ and $0<\theta<1$ or $q=1$ and $0\le\theta\le 1$\,. Then a function $f\in X_0+X_1$ belongs to $X_{\theta,q}=\big[X_0,X_1\big]_{\theta,q}$ if and only if $u=\sum_{k\ge 1}u_k$\,, where the series converges in $X_0+X_1$\,, and the sequence
\[
U_k=\mu_k^{-\vartheta}J(\mu_k,u_k)\in \ell_q(\NN)\,.
\]
In this case,  the norm $\|u\|_{X_{\theta,q},J}$ defined in \eqref{norm.interp.cont.1} is equivalent to
\begin{equation}\label{norm.interp.discr.1}
\|u\|_{\theta,q,JD}=\inf\Big\{\|U_k\|_{\ell^q(\NN)}\;\big|\; u=\sum_{k\ge 0}u_k\Big\}\,.
\end{equation}
\end{thm}
\noindent {\sl Proof.~}We prove the case $1\le q<\infty$ and we leave the easier case $q=\infty$ to the reader. We split the proof in two steps.

\noindent$\bullet~$\textsc{Step 1. }\textit{The inequality $\|u\|_{\theta,q,JD}\le c_1\|u\|_{X_{\theta,q},J}$}. First suppose that $u\in X_{\theta,q}$ and let $\varepsilon>0$\,. Then by definition \eqref{norm.interp.cont.1} of the norm $\|u\|_{X_{\theta,q},J}$,  there exists a function $f\in \LL^1_*$ such that
\begin{equation}\label{thm.JDiscr.step.1.1}
u=\int_0^{+\infty}f(t)\frac{\dt}{t}\qquad\mbox{and}\qquad \int_0^\infty\left(t^{-\theta}J(t,f(t))\right)^q\frac{\dt}{t}\le (1+\varepsilon)\|u\|_{X_{\theta,q},J}^q\,.
\end{equation}
Define the sequence $\{u_k\}_{k=1}^{\infty}$ by letting $\mu_0=0$ and
\[
u_k=\int_{\mu_k}^{\mu_{k+1}}f(t)\frac{\dt}{t}\qquad\mbox{so that}\qquad u=\sum_{k\ge 0}u_k
\]
the series being convergent in $X_0+X_1$\,, because the integral representation converges to $u$ there. Next, we observe that
\begin{equation}\label{thm.JDiscr.step.1.2}\begin{split}
J(\mu_k,u_k)
    &=\max\left\{\|u_k\|_{X_0}\,,\,\mu_k\,\|u_k\|_{X_1}\right\}
     \le \max\left\{\int_{\mu_k}^{\mu_{k+1}}\|f(t)\|_{X_0}\frac{\dt}{t}\,,\,
        \mu_k\,\int_{\mu_k}^{\mu_{k+1}}\|f(t)\|_{X_1}\frac{\dt}{t}\right\}\\
    &\le \int_{\mu_k}^{\mu_{k+1}}\max\left\{\|f(t)\|_{X_0}\,,\,
        \mu_k\,\|f(t)\|_{X_1}\right\}\frac{\dt}{t}
     \le \int_{\mu_k}^{\mu_{k+1}}\max\left\{\|f(t)\|_{X_0}\,,\,
        t\,\|f(t)\|_{X_1}\right\}\frac{\dt}{t}\\
     &=\int_{\mu_k}^{\mu_{k+1}}J(t,f(t))\frac{\dt}{t}\,,
\end{split}
\end{equation}
therefore,
\begin{equation}\label{thm.JDiscr.step.1.3}\begin{split}
\mu_k^{-\theta}\,J(\mu_k,u_k)
     &\le\int_{\mu_k}^{\mu_{k+1}}\mu_k^{-\theta}\, J(t,f(t))\frac{\dt}{t}
      =\frac{\mu_{k}^{-\theta}}{\mu_{k+1}^{-\theta}}
        \int_{\mu_k}^{\mu_{k+1}}\mu_{k+1}^{-\theta}\, J(t,f(t))\frac{\dt}{t}\\
     &\le_{(a)}\Lambda_0^\theta\int_{\mu_k}^{\mu_{k+1}}t^{-\theta}\, J(t,f(t))\frac{\dt}{t}\\
     &\le_{(b)}\Lambda_0^\theta\left(\int_{\mu_k}^{\mu_{k+1}}\frac{\dt}{t}\right)^{\frac{q-1}{q}}
        \left(\int_{\mu_k}^{\mu_{k+1}}\left(t^{-\theta}J(t,f(t))\right)^q\frac{\dt}{t}\right)^{\frac{1}{q}}\\
     &\le\Lambda_1 \left(\int_{\mu_k}^{\mu_{k+1}}\left(t^{-\theta}J(t,f(t))\right)^q\frac{\dt}{t}\right)^{\frac{1}{q}}
\end{split}
\end{equation}
where in $(a)$ we have used that $t\le \mu_{k+1}$ and that $\mu_{k+1}/\mu_{k}\le \Lambda_0<+\infty$\,. In $(b)$ we have used H\"older inequality, while in the last step we have used
\begin{equation}\label{thm.JDiscr.step.1.2b}
\left(\int_{\mu_k}^{\mu_{k+1}}\frac{\dt}{t}\right)^{\frac{q-1}{q}}
=\left(\log\frac{\mu_{k+1}}{\mu_k}\right)^{\frac{q-1}{q}}\le \left(\log\Lambda_0\right)^{\frac{q-1}{q}}
\end{equation}
and we have set $\Lambda_1=\Lambda_0^\theta\left(\log\Lambda_0\right)^{\frac{q-1}{q}} $\,. Summing on $k\ge 0$ inequality \eqref{thm.JDiscr.step.1.3} gives
\begin{equation}\label{thm.JDiscr.step.1.4}\begin{split}
\|U_k\|_{\ell^q(\NN)}&\le \sum_{k\ge 0}\left[\mu_k^{-\theta}\,J(\mu_k,u_k)\right]^q
      \le\Lambda_1 \sum_{k\ge 0}  \int_{\mu_k}^{\mu_{k+1}}\left(t^{-\theta}J(t,f(t))\right)^q\frac{\dt}{t} \\
     &=\Lambda_1   \int_{0}^{+\infty}\left(t^{-\theta}J(t,f(t))\right)^q\frac{\dt}{t}
     \le (1+\varepsilon)\|u\|_{X_{\theta,q},J}^q\\
\end{split}
\end{equation}
where in the last step we have used inequality \eqref{thm.JDiscr.step.1.1}\,. We obtain the desired inequality by letting $\varepsilon\to 0^+$\,.

\noindent$\bullet~$\textsc{Step 2. }\textit{The inequality $\|u\|_{\theta,q,J}\le c_0\|u\|_{X_{\theta,q},JD}$}. If $u=\sum_{k\ge 1}u_k$\,, where the series converges in $X_0+X_1$\,, we can define $u_0=0$ and a function $f\in\LL^1_*$ by
\begin{equation}\label{thm.JDiscr.step.2.1}
f(t)=\frac{u_k}{\log 2}\qquad\mbox{for}\qquad 2^k-1\le t\le 2^{k+1}-1\,,
\end{equation}
so that
\begin{equation}\label{thm.JDiscr.step.2.2}
\int_{2^k-1}^{2^{k+1}-1}f(t)\frac{\dt}{t}=u_k\qquad\mbox{and}\qquad u=\int_0^{+\infty}f(t)\frac{\dt}{t}\,.
\end{equation}
Notice that since $\mu_{k+1}/\mu_{k}\le \Lambda_0<+\infty$\,, we have
\begin{equation}\label{thm.JDiscr.step.2.3}\begin{split}
J(\mu_{k+1}\,,\,f(t))
&=\frac{1}{\log 2}\max\left\{\|u_k\|_{X_0}\,,\,\mu_{k+1}\,\|u_k\|_{X_1}\right\}\\
&\le\frac{\mu_{k+1}}{\mu_k\,\log 2}\max\left\{\|u_k\|_{X_0}\,,\,\mu_k\,\|u_k\|_{X_1}\right\}
\le \frac{\Lambda_0}{\log 2}J(\lambda_k,u_k)\,.
\end{split}
\end{equation}
Now consider the sequence $\mu_n$\,, and define the subsequences $\mu_{\underline{n_k}}$\,, and $\mu_{\overline{n_k}}$ through the indices
\[
\overline{n_k}=\min\left\{n\in \NN\;\big|\; \mu_n\ge 2^k-1\right\}\qquad\mbox{and}\qquad \underline{n_k}=\max\left\{n\in \NN\;\big|\; \mu_n\le 2^k-1\right\}
\]
Therefore, for all $k\ge 0$\,, we have that $\left[2^k,2^{k+1}\right]\subseteq\left[\lambda_{\underline{n_k}},\lambda_{\overline{n_k}}\right]$ and we can estimate:
\begin{equation}\label{thm.JDiscr.step.2.4}\begin{split}
\int_{2^k-1}^{2^{k+1}-1}\left(t^{-\theta}J(t,f(t))\right)^q\frac{\dt}{t}
&\le \sum_{n=\underline{n_k}}^{\overline{n_{k+1}}}\int_{\mu_n}^{\mu_{n+1}}\left(t^{-\theta}J(t,f(t))\right)^q\frac{\dt}{t}\\
&\le_{(a)} \sum_{n=\underline{n_k}}^{\overline{n_{k+1}}}
    \int_{\mu_n}^{\mu_{n+1}}\left(\mu_n^{-\theta}J(\mu_{n+1},f(t))\right)^q\frac{\dt}{t}\\
&\le_{(b)} \frac{\Lambda_0^q}{(\log 2)^q}
    \sum_{n=\underline{n_k}}^{\overline{n_{k+1}}}\left(\mu_n^{-\theta}J(\mu_n,u_n)\right)^q\int_{\mu_n}^{\mu_{n+1}}\frac{\dt}{t}\\
&\le_{(c)} \frac{\Lambda_0^q\log\Lambda_0}{(\log 2)^q}
    \sum_{n=\underline{n_k}}^{\overline{n_{k+1}}}\left(\mu_n^{-\theta}J(\mu_n,u_n)\right)^q
\end{split}\end{equation}
where in $(a)$ we have used that $t\ge\lambda_k$ and that $J$ is increasing with respect to the first variable $t\le \lambda_{k+1}$. In $(b)$ we have used inequality \eqref{thm.JDiscr.step.2.3}\,, while in the last step we have used inequality \eqref{thm.JDiscr.step.1.2b}\,.

Summing on $k\ge 0$ inequality \eqref{thm.JDiscr.step.2.4} gives
\begin{equation}\begin{split}
\int_{0}^{\infty}\left(t^{-\theta}J(t,f(t))\right)^q\frac{\dt}{t}
    &=\sum_{k\ge 0}\int_{2^k-1}^{2^{k+1}-1}\left(t^{-\theta}J(t,f(t))\right)^q\frac{\dt}{t}\\
      &\le \frac{\Lambda_0^q\log\Lambda_0}{(\log 2)^q}\sum_{k\ge 0}
      \sum_{n=\underline{n_k}}^{\overline{n_{k+1}}}\left(\mu_n^{-\theta}J(\lambda_k,u_k)\right)^q\\
     &\le_{(a)} 2\frac{\Lambda_0^q\log\Lambda_0}{(\log 2)^q}\sum_{k\ge 0}
        \left(\mu_n^{-\theta}J(\lambda_k,u_k)\right)^q = 2\frac{\Lambda_0^q\log\Lambda_0}{(\log 2)^q}\|U_k\|_{\ell^q(\NN)}^q\\
\end{split}
\end{equation}
where in $(a)$ we have used that
\[
\sum_{k\ge 0}
      \sum_{n=\underline{n_k}}^{\overline{n_{k+1}}}a_n
      + \sum_{k\ge 0}a_{\overline{n_k}}
      \le 2\sum_{k\ge 0}a_k
\]
since $a_{\overline{n_k}}$ is a subsequence of $a_k$\,. Summing up, we have found a function $f\in S(u)$ such that
\[
\int_{0}^{\infty}\left(t^{-\theta}J(t,f(t))\right)^q\frac{\dt}{t}\le 2\frac{\Lambda_0^q\log\Lambda_0}{(\log 2)^q}\|U_k\|_{\ell^q(\NN)}^q
\]
thus the above inequality holds for the norm $\|u\|_{\theta,q,J}$ defined in \eqref{norm.interp.cont.1}\,.\qed

\noindent\textbf{Remark. }There is another method that can be used to define the interpolation spaces, which is called the G-Method, see for instance \cite{Adams2003, LM}\,. Anyway, the two methods are equivalent, and the proof can be found in \cite{Adams2003} Thm. 7.16.

%%%%%%%%%%%%%%%%%%%%%%%%%%%%%%%%%%%%%%%%%%%%%%%%%%%%%%%%%%%%%%

\

\noindent {\large \sc Acknowledgment}

\noindent First and third author funded by Project MTM2011-24696 (Spain).  The second author would like to thank the hospitality of the Departamento de Matem\'{a}ticas of Universidad Aut\'{o}noma de Madrid. The second author is supported by the ANR project ``HAB''.

%%%%%%%%%%%%%%%%%%%%%%%%%%%%%%%%%%%%%%%%%%%%%%%%%%%%%%%%%%%%%%%%%%%%%

%%%%%%%%%%%%%%%%%%%%%%%%%%%%%%%%%%%%%%%%%%%%%%%%%%%%%%%%%%%%%%%%%%%%%

\bibliographystyle{amsplain}

\end{document}